\documentclass[]{article}

\usepackage{authblk}
\usepackage{amsmath}
\usepackage{amsfonts}
\usepackage{amsthm}
\usepackage{algpseudocode}
\usepackage{algorithm}
\usepackage{multirow}
\usepackage{hyperref}
\usepackage{cancel}
\usepackage{subcaption}
\usepackage{tikz} 
\usepackage{pgfplots}
\pgfplotsset{compat=1.16}
\graphicspath{{./Figures/}}

\newtheorem{definition}{Definition}

\usepackage{stmaryrd}

\newcommand{\domain}{\Omega}
\newcommand{\domainA}{\mathfrak{A}}
\newcommand{\domainB}{\mathfrak{B}}
\newcommand{\interface}{\mathfrak{I}}
\newcommand{\boundary}{\partial \Omega}
\newcommand{\bulk}{\Omega \setminus \mathfrak{I}}
\newcommand{\dens}{\rho}
\newcommand{\visc}{\mu}
\newcommand{\veloc}{\vectr{u}}
\newcommand{\press}{p}
\newcommand{\normal}{\vectr{n}}
\newcommand{\normalI}{\vectr{n}_{\mathfrak{I}}}
\newcommand{\normalGam}{\vectr{n}_{\Gamma}}

\newcommand{\normalB}{\vectr{n}_{\partial \Omega}}
\newcommand{\surfT}{\sigma}
\newcommand{\curv}{\kappa}
\newcommand{\levset}{\varphi}
\newcommand{\phiDG}{\varphi_{\textrm{evo}}}
\newcommand{\velocI}{\vectr{u}_{\mathfrak{I}}}
\newcommand{\grid}{\mathfrak{K}_{h}}
\newcommand{\cell}{K}
\newcommand{\grdSz}{h}
\newcommand{\cut}{\textrm{X}}
\newcommand{\edge}{\Gamma}
\newcommand{\edgeInt}{\Gamma_{\textrm{int}}}
\newcommand{\edgeD}{\Gamma_{\textrm{D}}}
\newcommand{\edgeN}{\Gamma_{\textrm{N}}}
\newcommand{\spc}{\mathfrak{s}}
\newcommand{\Pdeg}{k}
\newcommand{\brknPspacek}{\mathbb{P}_{k}}
\newcommand{\spaceVkX}{\mathbb{V}_{\vectr{k}}^{\textrm{X}}}
\newcommand{\basis}{\Phi}
\newcommand{\logEdg}{Edg}
\newcommand{\AggMap}{A}
\newcommand{\Agg}{Agg}
\newcommand{\AggTh}{\alpha}
\newcommand{\testV}{\vectr{v}}
\newcommand{\testP}{q}
\newcommand{\pnlty}{\eta}
\newcommand{\tangent}{\boldsymbol{\mathbf{\tau}}}
\newcommand{\projI}{\boldsymbol{\tensr{P}}_{\mathfrak{I}}}
\newcommand{\dt}{\Delta t}
\newcommand{\gridCC}{\mathfrak{K}_{\textrm{cc}}}
\newcommand{\gridNear}{\mathfrak{K}_{\textrm{near}}}
\newcommand{\gridNB}{\mathfrak{K}_{\textrm{nb}}}
\newcommand{\gridFar}{\mathfrak{K}_{\textrm{far}}}
\newcommand{\velocExt}{\vectr{u}_{\textrm{ext}}}
\newcommand{\velocEvo}{\vectr{u}_{\textrm{evo}}}
\newcommand{\massM}{\tensr{M}}
\newcommand{\OpM}{\tensr{Op}}
\newcommand{\RHS}{\textrm{RHS}}
\newcommand{\res}{\textrm{res}}
\newcommand{\ccNSE}{\epsilon_{\textrm{NSE}}}
\newcommand{\ccLS}{\epsilon_{\textrm{LS}}}
\newcommand{\laplace}{\textrm{La}}

\newcommand{\vectr}[1]{\boldsymbol{\mathbf{#1}}}
\newcommand{\tensr}[1]{\boldsymbol{\mathbf{#1}}}

\newcommand{\pDeriv}[2]{\frac{\partial{#1}}{\partial{#2}}}
\newcommand{\dif}[1]{~\textrm{d}{#1}}
\renewcommand{\div}[1]{\nabla \cdot {#1}}
\newcommand{\divI}[1]{\nabla_{\mathfrak{I}} \cdot {#1}}
\newcommand{\divH}[1]{\nabla_{h} \cdot {#1}}
\newcommand{\grad}[1]{\nabla{#1}}
\newcommand{\gradT}[1]{\nabla{#1}^{\textrm{T}}}
\newcommand{\gradH}[1]{\nabla_{h}{#1}}
\newcommand{\gradHT}[1]{\nabla_{h}{#1}^{\textrm{T}}}
\newcommand{\gradI}[1]{\nabla_{\interface}{#1}}
\newcommand{\jump}[1]{\left[\!\left[{{#1}}\right]\!\right]}
\newcommand{\aver}[1]{\left\{\!\left\{{#1}\right\}\!\right\}}
\newcommand{\norm}[1]{|\!| {#1} |\!|}
\newcommand{\LtwoNorm}[1]{|\!| {#1} |\!|_{L^2}}
\newcommand{\ltwoErr}[1]{|\!| {#1} |\!|_{l_2}}
\newcommand{\abs}[1]{| {#1} |}

\newcommand{\indA}[1]{{#1}_{\mathfrak{A}}}
\newcommand{\indB}[1]{{#1}_{\mathfrak{B}}}
\newcommand{\diri}[1]{{#1}_{\textrm{D}}}
\newcommand{\neum}[1]{{#1}_{\textrm{N}}}
\newcommand{\inn}[1]{{#1}^{-}}
\newcommand{\out}[1]{{#1}^{+}}
\newcommand{\slip}[1]{{#1}_{\textrm{slip}}}

\title{On a marching level-set method 
		for extended discontinuous Galerkin methods 
		for incompressible two-phase flows}
	
\author[1,2]{Martin Smuda}
\author[1,2]{Florian Kummer}

\affil[1]{Chair of Fluid Dynamics, TU Darmstadt, Germany}
\affil[2]{Graduate School of Computational Engineering, TU Darmstadt, Germany}

\begin{document}

\maketitle

\begin{abstract}
	In this work a solver for instationary two-phase flows on the basis of the extended Discontinuous Galerkin (extended DG/XDG) method is presented. The XDG method adapts the approximation space conformal to the position of the interface. This allows a sub-cell accurate representation of the incompressible Navier-Stokes equations in their sharp interface formulation.
The interface is described as the zero set of a signed-distance level-set function and discretized by a standard DG method. For the interface, resp. level-set, evolution an extension velocity field is used and a two-staged algorithm is presented for its construction on a narrow-band. On the cut-cells a monolithic elliptic extension velocity method is adapted and a fast-marching procedure on the neighboring cells. 
The spatial discretization is based on a symmetric interior penalty method and for the temporal discretization a moving interface approach is adapted. A cell agglomeration technique is utilized for handling small cut-cells and topology changes during the interface motion.
The method is validated against a wide range of typical two-phase surface tension driven flow phenomena including capillary waves, an oscillating droplet and the rising bubble benchmark.
\end{abstract}


\section{Introduction}


Within the past decades the significance of direct numerical simulations in context of multi-phase flows is increasing, not at least due to the increase of computational power. Up to date incompressible flow problems described by the Navier-Stokes equations (NSE) only allow analytical solutions for reduced or simplified settings. Considering in particular immiscible two-phase flows the complexity is even increased due to: discontinuous fluid properties across the interface, low-regularity solutions for flow properties and the presence of interfacial forces. Thus, sophisticated and specialized numerical methods are in need for the direct simulation of such problems.\\

Discontinuous Galerkin (DG) methods have gained quite some interest in the context of computational fluid dynamics in recent years. Reasons for that is the combination of favorable numerical properties known form other well established methods as the finite volume method (FVM) and the finite element method (FEM). DG discretizations are locally conservative due to the flux formulation and enables
a high-order method even on unstructured meshes by adapting the local ansatz functions. This cell locality is especially favorable for local mesh refinement and results in local mass matrices and sparse operator matrices. However, a major drawback of DG methods is the notably increase of degrees of freedom (DOF), when comparing to FVM on the same grid and interpolation order.

Considering the discretization of the sharp interface two-phase NSE with a DG method, the convergence order degenerates for such low-regularity solutions (i.e. kinks on the velocity field and jumps in the pressure field) if the approximation space is not conformal to the interface position. However, in context of transient flow problems with moving interfaces one wants to avoid costly remeshing every time step in order to adapt the solution space. A first method that enriches the approximation space such that the solution may exhibit discontinuities inside a finite element is presented in Moës et al.\cite{moes_finite_1999} introducing the extended FEM (XFEM) for the simulation of crack growth in solid mechanics. The first application to incompressible two-phase flows is found in Groß and Reusken\cite{gros_extended_2007}, where the pressure is discretized with XFEM and exhibits a jump due to the surface tension force. A discretization of both pressure and velocity with XFEM is provided by Fries\cite{fries_intrinsic_2009}, which additionally allows the representation of the kink in the velocity field.

The first extended method for DG is presented in Bastian and Engwer\cite{bastian_unfitted_2009} for the discretization of elliptic scalar problems. In Heimann et al.\cite{heimann_unfitted_2013} this approach is applied to incompressible Navier-Stokes two-phase flows. The unfitted DG (UDG) method is based on the non-symmetric interior penalty method, where the cut out meshes for both subdomains are based on a piece-wise linear approximation of the interface. The XDG method proposed by Kummer\cite{kummer_extended_2016} for steady two-phase flows uses a high-order approximation of the interface in combination with a quadrature technique (i.e. the Hierarchical Moment Fitting\cite{muller_highly_2013}) for the implicitly defined domains. The discretization is based on the symmetric interior penalty (SIP) method\cite{arnold_interior_1982} and the stabilization against small cut-cells is ensured by cell agglomeration. Another approach by Saye\cite{saye_implicit_2017} is the use of implicitly defined meshes with curved elements that are interface-conforming. The implicit mesh DG method provides high-order accuracy for the interface jump conditions in combination with the interfacial gauge method\cite{saye_interfacial_2016}, where the numerical coupling between the fluid velocity, pressure, and the interface position is reduced.\\

In this work we follow the XDG method according to Kummer\cite{kummer_extended_2016}. For the extension to the transient two-phase problem, we apply the XDG adapted moving interface approach introduced by Kummer et al.\cite{kummer_time_2018}. This temporal discretization allows for high orders of accuracy in time without a full discretization on space-time elements as done e.g. in Lehrenfeld\cite{lehrenfeld_nitsche_2015} and in Cheng and Fries\cite{cheng_higher-order_2009} for XFEM. The moving interface approach allows to keep the spatial discretization unchanged and to obtain a conservative scheme in time. The conservation property is not given for splitting approaches such as the Strang splitting used in Heimann et al.\cite{heimann_unfitted_2013}. In contrast to Kummer\cite{kummer_extended_2016} the quadrature method proposed by Saye\cite{saye_high-order_2015} is used. Furthermore, the numerical treatment of the surface tension force is in this work handled via the Laplace-Beltrami formulation.  

For the representation of the interface we employ a signed-distance level-set field. The numerical discretization is done by a standard DG formulation resulting in a high-order and sub-cell accurate approximation of the interface. For the evolution of the level-set field a extension velocity field is used in order to preserve the signed-distance property and reduce the need for reinitialization. 
For the construction of such an extension velocity field, one can employ direct approaches\cite{chessa_extended_2002, gibou_fourth_2005, sauerland_stable_2013}, fast-marching\cite{adalsteinsson_fast_1999}/fast-sweeping methods\cite{aslam_static_2014} or PDE-based techniques\cite{chen_simple_1997, utz_high-order_2018}. In this work we present a two-staged construction with an high-order elliptic PDE-based approach on the interface cells and a low-order fast-marching on the neighboring cells. 

The range of level-set reinitialization methods is similar to the one for the extension velocity problem. It includes direct geometry-based methods\cite{marchandise_stabilized_2007, saye_high-order_2015}, iterative fast-marching\cite{adalsteinsson_fast_1995}/fast-sweeping\cite{zhao_fast_2004} methods and PDE-based approaches\cite{sussman_level_1994, chunming_li_level_2005, basting_minimization-based_2013}. For our method we apply again an elliptic PDE-based reinitialization\cite{utz_interface-preserving_2017} on the interface cells and perform a fast-marching procedure on the neighboring cells.\\

The outline of this work is as follows. In Section \ref{sec:problem} the continuous problem statement for the transient incompressible two-phase flow for immiscible fluids given. An introduction to the XDG space is provided in Section \ref{sec:XDGmain}. The numerical representation of the level-set field and its evolution by the two-staged evolution algorithm is considered in Section \ref{sec:LevelSet}. Here, we differentiate between the
interface evolution in context of the XDG method (Section \ref{sec:InterfaceEvolution}) and the fast-marching algorithm for the level-set field (Section \ref{sec:FastMarchingEvo}). The XDG discretization of the considered problem is given in Section \ref{sec:XDGmethod}. Furthermore, we briefly outline the implemented cell-agglomeration technique in terms of graph theory. The over-all solver structure and the adaptive local mesh refinement are described in Section \ref{sec:Solver}. The presented solver is validated against a range of typical two-phase surface tension driven flows in Section \ref{sec:results}. We conclude this work in Section \ref{sec:conclusion}.

%
\section{Problem statement - Transient incompressible two-phase flow for immiscible fluids with a dividing sharp interface}
\label{sec:problem}

We consider the two-phase setting within a sharp interface formulation and define the computational domain $\domain \subset \mathbb{R}^2$ as the disjoint partitioning of the time-dependent fluid bulk phases $\domainA(t)$ and $\domainB(t)$ and the moving interface $\interface(t)$ by
\begin{align}
	\domain = \domainA(t) \mathbin{\dot{\cup}} \interface(t) \mathbin{\dot{\cup}} \domainB(t).
	\label{eq:twoPhaseSetting}
\end{align}
We restrict ourselves to the two-dimensional case and assume that the sharp interface  $\interface(t) := \overline{\domainA}(t) \cap \overline{\domainB}(t)$ is a one-dimensional manifold. The considered two-phase flow problems are described by the transient incompressible Navier-Stokes equations (NSE) denoting the momentum balance equation \eqref{eq:momentum} and the continuity equation \eqref{eq:continuity} by
\begin{subequations}
\begin{align}
	\dens \pDeriv{\veloc}{t} + \dens \div{\left(\veloc \otimes \veloc \right)} &= \div{\left( -\press \tensr{I} + \visc \left( \grad{\veloc} + \gradT{\veloc} \right) \right)} + \dens \vectr{g} \quad \text{   in } \bulk, \label{eq:momentum}\\
	\div{\veloc} &= 0 \quad \quad \quad \quad \quad \quad \quad \quad \quad \quad \quad \quad \quad \quad \text{in } \bulk, \label{eq:continuity}
\end{align}
\label{eq:NavierStokes}
\end{subequations}
where $\veloc = \veloc(\vectr{x},t)$ is the velocity vector, $\press = \press(\vectr{x},t)$ the pressure and $\vectr{g} = \vectr{g}(\vectr{x},t)$ describes a volume force, e.g. gravity. The physical properties density $\dens = \dens(\vectr{x})$ and dynamic viscosity $\visc = \visc(\vectr{x})$ are piece-wise constant defined by
\begin{equation}
    \dens(\vectr{x}) = 
    \begin{cases} 
    \indA{\dens} \text{ for } \vectr{x} \in \domainA\\ 
    \indB{\dens} \text{ for } \vectr{x} \in \domainB 
    \end{cases} 
    \text{and} \quad
    \visc(\vectr{x}) = 
    \begin{cases} 
    \indA{\visc} \text{ for } \vectr{x} \in \domainA\\ 
    \indB{\visc} \text{ for } \vectr{x} \in \domainB 
    \end{cases}. 
\end{equation}
The corresponding jump conditions for the NSE \eqref{eq:NavierStokes} at a material interface $\interface$ are given by
\begin{subequations}
\begin{align}
	\jump{ -\press \normalI + \visc \left( \grad{\veloc} + \gradT{\veloc} \right) \normalI} &= \surfT \curv \normalI \quad \text{    on } \interface, \label{eq:jump_momentumMaterial}\\
	\jump{\veloc} &= 0 \quad \quad \quad \text{ on } \interface. \label{eq:jump_massMaterial}
\end{align}
\label{eq:jump_NavierStokes}
\end{subequations}
The jump operator $\jump{\cdot}$ is defined in \eqref{eq:JumpOp} and the interface normal $\normalI$ is pointing from $\domainA$ to $\domainB$. In the momentum jump condition \eqref{eq:jump_momentumMaterial} on the right-hand side $\surfT$ denotes the surface tension coefficient and $\curv$ the mean curvature of the interface $\interface$. At the domain boundary $\boundary = \diri{\boundary} \mathbin{\dot{\cup}} \neum{\boundary}$, describing a disjoint decomposition into Dirichlet $\diri{\boundary}$ and Neumann regions $\neum{\boundary}$, the corresponding boundary conditions read
\begin{subequations}
	\begin{align}
	\veloc &= \diri{\veloc} \quad \text{    on } \diri{\boundary}, \label{eq:Dirichlet}\\
	-\press \normal_{\boundary} + \visc \left( \grad{\veloc} + \gradT{\veloc} \right) \normal_{\boundary} &= 0 \quad \quad \text{ on } \neum{\boundary} \label{eq:Neumann}
	\end{align}
	\label{eq:boundary_conditions}
\end{subequations}
with a given function $\diri{\veloc} = \diri{\veloc}(\vectr{x}, t)$ on $\diri{\boundary}$. In order to close the initial value problem we set an initial condition for the velocity field by 
\begin{equation}
    \veloc(\vectr{x},0) = \veloc_0(\vectr{x}) \quad \textrm{ with } \div{\veloc_0} = 0 \quad \textrm{ for } \vectr{x} \in \domain \setminus \interface(0),
\end{equation}
where the initial interface position $\interface(0)$ is given. The material interface evolves according to the bulk velocity $\veloc(\vectr{x},t)$ at $\vectr{x} \in \interface(t)$.

%
\section{The extended Discontinuous Galerkin space}
\label{sec:XDGmain}

In this section the extended Discontinuous Galerkin (extended DG/XDG) space according to Kummer \cite{kummer_extended_2016} is briefly introduced. In this case the interface is represented by a sufficiently smooth level-set function $\levset(\vectr{x}, t)$ that is almost everywhere $C^1(\domain)$-continuous, but at least $\levset(\vectr{x}, t)\in C^0(\domain)$. Thus, the partitioning of the computational domain $\domain$ in \eqref{eq:twoPhaseSetting} is implicitly defined by
\begin{subequations}
	\begin{align}
	\domainA(t) &= \{ \vectr x \in \domain : \levset(\vectr x, t) < 0 \},\\
	\domainB(t) &= \{ \vectr x \in \domain : \levset(\vectr x, t) > 0 \},\\
	\interface(t) &= \{ \vectr x \in \domain : \levset(\vectr x, t) = 0 \}.
	\end{align}
	\label{eq:defLevSet}
\end{subequations}
Such a representation of the interface via a level-set function allows the direct computation of the interface normal $\normalI$ by
\begin{equation}
	\normalI = \frac{\grad{\levset}}{\abs{\grad{\levset}}}.
	\label{eq:LSnormal}
\end{equation}
The numerical discretization of $\levset$ and its evolution in context of the XDG method are discussed in the next Section \ref{sec:LevelSet}. In order to formulate the XDG space, first, some basic definitions and operators are introduced. We assume that the computational domain $\domain$ is polygonal and simply connected. Thus, the following basic definitions framework hold.
\begin{definition}[basic notations] For a polygonal and simply connected computational domain $\domain = \domainA(t) \mathbin{\dot{\cup}} \interface(t) \mathbin{\dot{\cup}} \domainB(t) \in \mathbb{R}^2$ we define:
	\begin{itemize}
		\item{the numerical background mesh $\grid = \{ \cell_1, ... ,\cell_J\}$ that covers the whole domain $\overline{\domain} = \cup_j \overline{\cell}_j$ with non-overlapping cells ($\int_{\cell_j\cap\cell_l} 1 \dif{V} = 0, l \neq j$), where $\grdSz$ denotes the maximum diameter of all cells $\cell_j$. 
		}
		\item{the set of all edges in the mesh is given by $\edge(t) := \cup_j \partial \cell_j \cup \interface(t)$. This set can be subdivided into $\edge(t) = \edgeInt(t) \cup \edgeD \cup \edgeN$, where $\edgeInt(t) = \edge(t) \setminus \boundary$ denotes the set of all internal edges, $\edgeD$ and $\edgeN$ the set of edges imposed with Dirichlet and Neumann boundary conditions, respectively.
		}
		\item{a normal field $\normalGam$ on $\edge$, where it is $\normalGam = \normalI$ on $\interface$ and represents an outer normal on $\boundary$ with $\normalGam = \normal_{\boundary}$. 
		}
		\item{the broken gradient $\gradH{\vectr{f}}$ defines for $\vectr{f} \in C^1(\domain \setminus \edge)$ the gradient on the domain $\domain \setminus \edge$. According to that, the broken divergence $\divH{\vectr{f}}$ is defined.
		}
	\end{itemize}
\label{def:baiscNotations}
\end{definition}
The introduction of the interface $\interface$ within the two-phase setting allows the definition of cut-cells and a corresponding cut-cell mesh.
\begin{definition}[cut-cell and cut-cell mesh] For a numerical background mesh $\grid = \{ \cell_1, ..., \cell_J\}$ we define the time-dependent cut-cells as
\begin{equation}
    \cell^{\cut}_{j,\spc}(t) := \cell_j \cap \spc(t)
\end{equation}
for a species $\spc(t) \in \{ \domainA(t), \domainB(t)\}$. The set of of all cut-cells defines the time-dependent cut-cell mesh
\begin{equation}
    \grid^{\cut}(t) = \{ \cell^{\cut}_{1,\domainA}(t), \cell^{\cut}_{1,\domainB}(t), ..., \cell^{\cut}_{J,\domainA}(t), \cell^{\cut}_{J,\domainB}(t) \}.
\end{equation}
Note that in a background cell $\cell_j$ where the interface is located we end up with two cut-cells $\cell_{j,\domainA}^{\cut}$ and $\cell_{j,\domainB}^{\cut}$, otherwise the original background cell is recovered for the respective species $\spc$, i.e. $\cell^{\cut}_{j,\domainA}(t) = \cell_j$ and $\cell^{\cut}_{j,\domainB}(t) = \emptyset$, or $\cell^{\cut}_{j,\domainA}(t) = \emptyset$ and $\cell^{\cut}_{j,\domainB}(t) = \cell_j$.
\end{definition}
The XDG method is essentially a DG method applied on the cut-cell mesh. Thus, we define the corresponding XDG space as follows.
\begin{definition}[XDG space] The broken cut-polynomial space $\brknPspacek^{\cut}$ of total degree $\Pdeg$ is defined as
	\begin{equation}
	\begin{aligned}
	\brknPspacek^{\cut}(\grid, t) := \left\{ f \in L^2(\domain); \forall \cell \in \grid : f\vert_{\cell \cap \spc(t)} \text{ are polynomial, }\right.\\ \left. 
	\text{deg}(f\vert_{\cell \cap \spc(t)}) \leq \Pdeg \right\} = \brknPspacek(\grid^{\cut}(t)),  
	\end{aligned}
	\end{equation}
	where $\brknPspacek$ denotes the standard broken polynomial space of total degree $\Pdeg$.
\end{definition}
The corresponding local approximation of a field property $f_j$ in a background cell $\cell_j$ occupied by both phases is given by 
\begin{equation}
		f_j(\vectr{x},t) = \sum_{l=0}^{k} \tilde{f}_{j,l,\domainA}(t) \basis_{j,l}(\vectr{x})\textrm{1}_{\domainA}(\vectr{x},t) + \tilde{f}_{j,l,\domainB}(t) \basis_{j,l}(\vectr{x})\textrm{1}_{\domainB}(\vectr{x},t), \quad \vectr{x} \in \cell_j,
	\label{eq:localApproxCC}
\end{equation}
where $\textrm{1}_{\spc}(\vectr{x},t)$ defines the characteristic function for species $\spc(t)$, i.e. $\textrm{1}_{\spc}(\vectr{x},t) = 1$ for $\vectr{x} \in \spc(t)$ and zero everywhere else. Thus, the corresponding cut-polynomial basis functions for $\spc$ are given by $ \basis_{j,l}(\vectr{x}) \textrm{1}_{\spc}(\vectr{x},t) = \basis_{j,l,\spc}^{\cut}(\vectr{x},t) \in \brknPspacek^{\cut}(\grid, t)$. Note that the time-independent local basis functions $\basis_{j,l} \in \brknPspacek(\grid)$ are not altered. In such cut-cells we provide for each phase separate coefficients $\tilde{f}_{j,l,\domainA}(t)$ and $\tilde{f}_{j,l,\domainB}(t)$, also denoted as degrees of freedom (DOF). 

In order to formulate an XDG discretization for the considered problem (Section \ref{sec:XDGdiscretization}) we introduce at last the jump and average operator on the edges $\edge$.
\begin{definition}[inner and outer value, jump and average operator] For a field $f \in \mathcal{C}^0(\domain \setminus \edge)$ we define the inner and outer values, $\inn{f}$ and $\out{f}$, at the edges $\edge$ by:
	\begin{subequations}
		\begin{align}
		\inn{f}(\vectr{x}) &:= \lim\limits_{\xi \searrow 0} f(\vectr{x} - \xi \normalGam) \quad \text{ for } \vectr{x} \in \edge,
		\label{eq:uin}\\
		\out{f}(\vectr{x}) &:= \lim\limits_{\xi \searrow 0} f(\vectr{x} + \xi \normalGam) \quad \text{ for } \vectr{x} \in \edgeInt.
		\label{eq:uout}
		\end{align}
	\end{subequations}
	Thus, the jump and average operator are defined as: 
	\begin{eqnarray}
	\jump{f} & := &
	\left\{  \begin{array}{ll} 
	\inn{f} - \out{f} &  \text{on } \edgeInt       \\
	\inn{f}            &  \text{on } \boundary 
	\end{array} \right. ,
	\label{eq:JumpOp}
	\\
	\aver{f} & := &
	\left\{  \begin{array}{ll} 
	\frac{1}{2}(\inn{f}+ \out{f})  &  \text{on } \edgeInt       \\
	\inn{f}                &  \text{on } \boundary
	\end{array} \right.  .
	\label{eq:AverOP}
	\end{eqnarray}
\end{definition}

%
\section{Level-set discretization and interface evolution}
\label{sec:LevelSet}


We define the level-set function by its signed-distance formulation leading to the property that $\abs{\grad\levset} = 1$. In order to allow a high-order and sub-cell accurate interface representation, the level-set function is approximated by a standard DG space $\mathbb{P}_k(\grid)$. However, such a discretization is in general discontinuous and not suitable for a spatial discretization in the XDG space, which requires at least that $\levset(\vectr{x}, t)\in C^0(\domain)$. Thus, we introduce two level-set fields 
\begin{equation}
    \phiDG \in \mathbb{P}_k(\grid) \quad \textrm{and} \quad \levset \in \mathbb{P}_{k+1}(\grid).
    \label{eq:LevelSets}
\end{equation}
The first one is used for handling the level-set evolution (Section \ref{sec:FastMarchingEvo}) and we choose the same polynomial degree as for the velocity field with $\Pdeg$ (see Eq. \eqref{eq:functionSpace}). The second level-set field is used for the discretization (e.g. the interface normal $\normalI$ and numerical integration on cut-cells) and is given as the $L^2$-projection of $\phiDG$ constrained by continuity conditions at the cell boundaries (Section \ref{sec:InterfaceEvolution}). The polynomial degree is chosen to be $\Pdeg + 1$.

\begin{figure}
	\centering
	\def\svgwidth{240pt}
\begingroup%
  \makeatletter%
  \providecommand\color[2][]{%
    \errmessage{(Inkscape) Color is used for the text in Inkscape, but the package 'color.sty' is not loaded}%
    \renewcommand\color[2][]{}%
  }%
  \providecommand\transparent[1]{%
    \errmessage{(Inkscape) Transparency is used (non-zero) for the text in Inkscape, but the package 'transparent.sty' is not loaded}%
    \renewcommand\transparent[1]{}%
  }%
  \providecommand\rotatebox[2]{#2}%
  \newcommand*\fsize{\dimexpr\f@size pt\relax}%
  \newcommand*\lineheight[1]{\fontsize{\fsize}{#1\fsize}\selectfont}%
  \ifx\svgwidth\undefined%
    \setlength{\unitlength}{441.55977944bp}%
    \ifx\svgscale\undefined%
      \relax%
    \else%
      \setlength{\unitlength}{\unitlength * \real{\svgscale}}%
    \fi%
  \else%
    \setlength{\unitlength}{\svgwidth}%
  \fi%
  \global\let\svgwidth\undefined%
  \global\let\svgscale\undefined%
  \makeatother%
  \begin{picture}(1,0.37808314)%
    \lineheight{1}%
    \setlength\tabcolsep{0pt}%
    \put(0,0){\includegraphics[width=\unitlength,page=1]{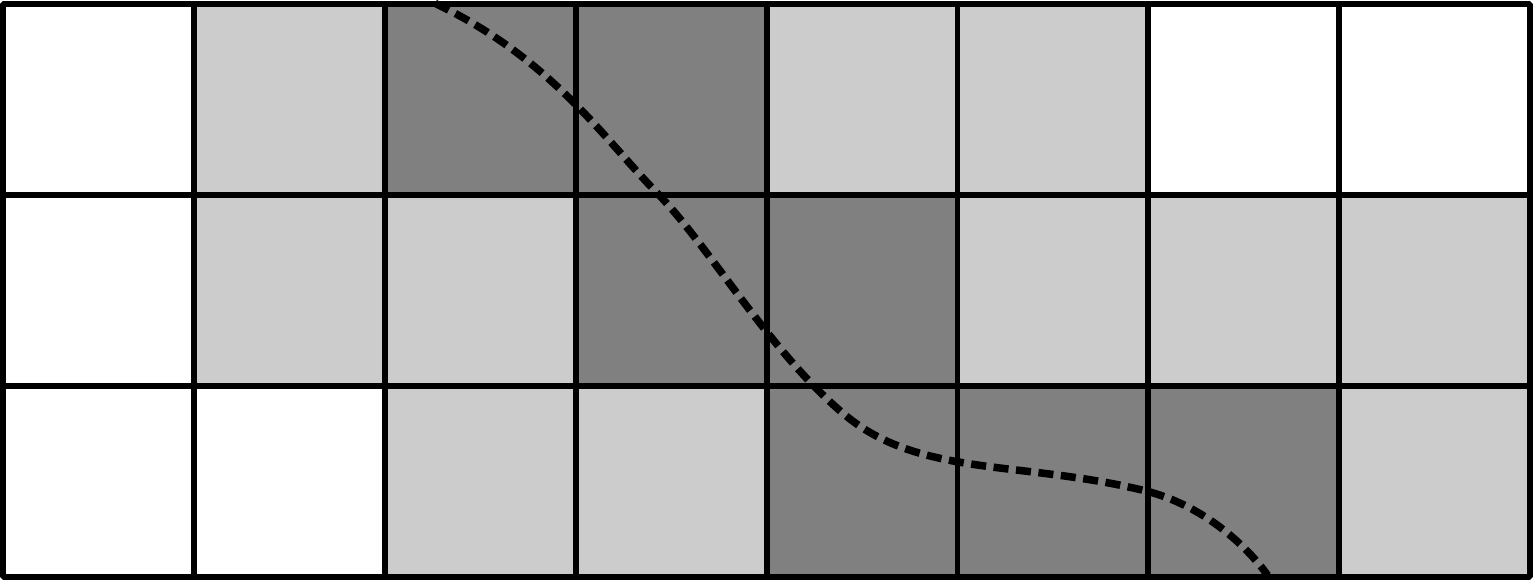}}%
    \put(0.81206564,0.05875264){\color[rgb]{0,0,0}\makebox(0,0)[t]{\lineheight{1.25}\smash{\begin{tabular}[t]{c}$\interface$\end{tabular}}}}%
    \put(0.56372445,0.18213617){\color[rgb]{0,0,0}\makebox(0,0)[t]{\lineheight{1.25}\smash{\begin{tabular}[t]{c}$\gridCC$\end{tabular}}}}%
    \put(0.31193796,0.17991013){\color[rgb]{0,0,0}\makebox(0,0)[t]{\lineheight{1.25}\smash{\begin{tabular}[t]{c}$\gridNear$\end{tabular}}}}%
    \put(0.06454436,0.1831745){\color[rgb]{0,0,0}\makebox(0,0)[t]{\lineheight{1.25}\smash{\begin{tabular}[t]{c}$\gridFar$\end{tabular}}}}%
  \end{picture}%
\endgroup%

	\caption{The narrow-band $\gridNB := \gridCC \cup \gridNear$ around the interface includes cut-cells $\gridCC$ (dark gray cells) and the neighbors $\gridNear$ (light gray) considering the vertices of the cut-cells. The white cells are denoted as far-field cells $\gridFar$.}
	\label{fig:narrowband}
\end{figure}

The following operations related to the level-set, resp. interface, evolution are only performed on a subset around the interface denoted as the narrow-band $\gridNB$. This narrow-band $\gridNB := \gridCC \cup \gridNear$ includes the set of cut-cells $\gridCC$ and its point-neighbors $\gridNear$, see Figure \ref{fig:narrowband}. According to that, both level-set fields $\phiDG$ and $\levset$ exhibits the signed-distance property only on the narrow-band. On the far field cells $\gridFar$ the value is either set to $\phiDG = -1$ for cells in $\domainA$ or $\phiDG = +1$ in $\domainB$.

\subsection{Interface evolution algorithm for XDG - Enforcing continuity of the level-set field}
\label{sec:InterfaceEvolution}

The evolution equation for the level-set field $\phiDG$ in the computational domain $\domain$ is given by 
\begin{equation}
    \pDeriv{\phiDG}{t} + \velocEvo \cdot \grad{\phiDG} = 0 \quad \textrm{ in } \domain,
    \label{eq:levsetEvo}
\end{equation}
where the evolution velocity for a material interface $\interface$ is set to $\velocEvo = \veloc$. However, in context of the XDG method the approximation space for the bulk velocity is in $\brknPspacek^{\cut}(\grid)$, whereas the level-set evolution needs a velocity field $\velocEvo \in \brknPspacek(\grid)$ according to $\phiDG$. In order to generate a velocity field for the advection equation \eqref{eq:levsetEvo}, we introduce a density-averaged velocity field $\veloc_{\dens \textrm{Aver}} \in \brknPspacek(\grid)$, where the velocity components $u_j^{d}$, $d=\{1,2\}$, in each cell $\cell_j \in \gridNB$ are constructed as   \begin{equation}
	u_j^{d}(\vectr{x}) = \frac{1}{\indA{\dens} + \indB{\dens}} \sum_{l=0}^{k} \left( \indA{\dens} \tilde{u}_{j,l,\domainA}^{d} \basis_{j,l}(\vectr{x}) + \indB{\dens} \tilde{u}_{j,l,\domainB}^{d} \basis_{j,l}(\vectr{x}) \right), \quad \vectr{x} \in \cell_j, \basis_{j,l} \in \brknPspacek(\grid).
	\label{eq:meanVelocity}
\end{equation}
The coefficients $\tilde{u}_{j,l,\domainA}^{d}$ and $\tilde{u}_{j,l,\domainB}^{d}$ are the same coefficients for $\veloc \in \brknPspacek^{\cut}$. Recalling the local XDG approximation in \eqref{eq:localApproxCC}, we omit the identity function to recover the standard approximation in $\brknPspacek(\grid)$. A marching procedure is employed for the level-set advection and is presented in Section \ref{sec:FastMarchingEvo}.

\paragraph{Enforcing continuity of the level-set field $\levset$}

Since we are dealing with two level-set fields $\phiDG$ and $\levset$, the latter one needs to be constructed from $\phiDG$ after the advection. 
Therefore, we employ a $L^2$-projection of $\phiDG$ onto $\levset$ with additional continuity constraints that is described as the following quadratic optimization problem (OP) with 
\begin{subequations}
	\begin{align}
	&\min \LtwoNorm{\levset - \phiDG}^2 \quad \textrm{ on } \gridNB\\
	&\text{s.t.} \quad \inn{\levset}\vert_{\edge_i} = \out{\levset}\vert_{\edge_i}, \quad \forall \edge_i \in \edge_{\textrm{int}, \textrm{nb}}.
	\end{align}
	\label{eq:L2ContinuityProjection}
\end{subequations}
The projection is restricted to the narrow-band $\gridNB$ and the continuity conditions are defined on the set of internal edges in the narrow-band $\edge_{\textrm{int}, \textrm{nb}}$. Let be $\levset = \tilde{\vectr{\levset}}^{\textrm{T}}\vectr{\basis}$ with $ \vectr{\basis} \in \mathbb{P}_{k+1}(\grid)$ and $\tilde{\vectr{\levset}}$ the sought-after coordinate vector and $\phiDG^{\ast} = \tilde{\vectr{\varphi}}_{\textrm{evo}}^{\textrm{T}}\vectr{\basis}$ the projection of $\phiDG \in \mathbb{P}_{k}(\grid)$ onto $\mathbb{P}_{k+1}(\grid)$. Thus, the OP is equivalently stated as 
\begin{subequations}
\begin{align}
    &J(\tilde{\vectr{\levset}}) = \frac{1}{2} \tilde{\vectr{\levset}}^{\textrm{T}} \tensr{M} \tilde{\vectr{\levset}} - \tilde{\vectr{\levset}}^{\textrm{T}} \tilde{\vectr{\varphi}}_{\textrm{evo}} + d \rightarrow \textrm{ min },\\
    &\text{s.t.} \quad\tensr{A} \tilde{\vectr{\levset}} = 0 \label{eq:equalityConstraint}
\end{align}
\label{eq:OPmatrixForm}
\end{subequations}
where $\tensr{M} = \tensr{I}$, since $\vectr{\basis}$ is a orthonormal basis, $d = \tilde{\vectr{\varphi}}_{\textrm{evo}}^{\textrm{T}} \tilde{\vectr{\varphi}}_{\textrm{evo}}$ and $\tensr{A}$ describes the constraint matrix. The equality constraints \eqref{eq:equalityConstraint} are incorporated via Lagrange multipliers $\vectr{\lambda}$ following the Karush-Kuhn-Tucker conditions of first order. The solution of the resulting OP can be described as the solution of the following system of linear equations
\begin{subequations}
\begin{align}
    \tensr{M}\tilde{\vectr{\levset}} 
    - \tilde{\vectr{\varphi}}_{\textrm{evo}}
    + \tensr{A}^{\textrm{T}} \vectr{\lambda} &= 0\\
    \tensr{A}\tilde{\vectr{\levset}} &= 0.
\end{align}
\end{subequations}
This saddle point system is solved via the Schur complement method, The Schur complement in this case is given by $\tensr{A} \tensr{M}^{-1} \tensr{A}^{\textrm{T}} = \tensr{A} \tensr{A}^{\textrm{T}}$ and describes a symmetric positive definite matrix. Thus, the saddle point system is reduced to solve 
\begin{equation}
    \tensr{A} \tensr{A}^{\textrm{T}} \vectr{\lambda} = \tensr{A} \tilde{\vectr{\varphi}}_{\textrm{evo}}
    \label{eq:OPschur}
\end{equation}
for $\vectr{\lambda}$ and evaluate the solution vector $\tilde{\vectr{\levset}}$ via
\begin{equation}
    \tilde{\vectr{\levset}} = \tilde{\vectr{\varphi}}_{\textrm{evo}} 
    - \tensr{A}^{\textrm{T}} \vectr{\lambda}.  
    \label{eq:OPsolution}
\end{equation}
Note that the problem size of \eqref{eq:OPschur} is reduced to the number of the constraint conditions. The continuity constraints \eqref{eq:equalityConstraint} at the internal edges $\edge_{\textrm{int}, \textrm{nb}}$ are given as the following equality conditions at a sufficient amount of points $pt$ at $\edge_i$ with
\begin{equation}
    \sum_{n=0}^{k+1} \tilde{\levset}_{i,n}^{-} \basis_{i,n}^{-}(pt) - \sum_{n=0}^{k+1} \tilde{\levset}_{i,n}^{+} \basis_{i,n}^{+}(pt) = 0,
    \label{eq:ConstraintMatrix}
\end{equation}
where $\tilde{\levset}_{i,n}^{-}$ and $\tilde{\levset}_{i,n}^{+}$ denote the coefficients of the corresponding inner and outer cell at $\edge_i$. The condition \eqref{eq:ConstraintMatrix} represents one row in the constraint matrix A. The amount of conditions is chosen according to the approximation space of $\levset = \mathbb{P}_{k+1}(\grid)$. Such a formulation of suitable equality conditions for the continuity projection is also applicable on meshes with hanging nodes, which is necessary for the implementation of adaptive mesh refinement as described in Section \ref{sec:Solver}.

\paragraph{Interface evolution algorithm for a XDG method}

The entire interface evolution algorithm including the update of the integration domains and the narrow-band is given in Algorithm \ref{alg:UpdateLevelSet}. The actual advection of the level-set field $\phiDG$ is performed in \textproc{MarchingLevelSetEvolution} (line \ref{algLine:marchingLevSetEvo}) and is presented in the next section. 
\begin{algorithm}
	\caption{Interface evolution algorithm}
	\label{alg:UpdateLevelSet}
	\begin{algorithmic}[1]
		\Procedure{UpdateLevelSet}{$\veloc^{n}$, $\dt$}
		\State $\gridCC^{n}, \gridNear^{n}, \gridFar^{n} \gets \gridCC^{n+1}, \gridNear^{n+1}, \gridFar^{n+1}$
		\State $\veloc_{\dens \textrm{Aver}} \gets \textrm{GetAveragedVelocityFromXDGField}(\veloc^{n}, \gridNB^{n})$	\Comment{density-averaged velocity \eqref{eq:meanVelocity}}
		\State $\phiDG^{n+1} \gets \textproc{MarchingLevelSetEvolution}(\phiDG^{n}, \label{algLine:LevelSetEvolution} \veloc_{\dens \textrm{Aver}}, \dt)$ \Comment{Algorithm \ref{alg:marchingLevSetEvo}} \label{algLine:marchingLevSetEvo}
		\State $\levset^{n+1} \gets \textrm{ContinuityProjection}(\phiDG^{n+1}, \gridNB^{n})$ \Comment{Solve Eq. \eqref{eq:OPschur} and evaluate Eq. \eqref{eq:OPsolution}}
		\State $\domainA^{n+1}, \domainB^{n+1}, \interface^{n+1} \gets \textrm{UpdateDomains}(\levset^{n+1})$ 
		\State $\gridCC^{n+1}, \gridNear^{n+1}, \gridFar^{n+1} \gets \textrm{UpdateNarrowBand}(\levset^{n+1})$ 
		\State \textbf{return} $\levset^{n+1}$
		\EndProcedure
	\end{algorithmic}
\end{algorithm}

\subsection{A marching level-set advection method}
\label{sec:FastMarchingEvo}

Considering the advection of the level-set field $\phiDG$ by equation \eqref{eq:levsetEvo} with $\velocEvo = \veloc$, this would not preserve the desired signed-distance property $\abs{\grad\phiDG} = 1$. Furthermore, one should note that for the advection of the interface the only valid velocity for the evolution should be the one at the interface itself. Thus, in order to preserve the desired signed distance property, at least sufficiently far from non-smooth regions in $\phiDG$ and for some finite time, we construct the evolution velocity $\velocEvo$ as an extension velocity field $\velocExt$. Therefore, we solve the extension velocity problem (EVP) for each velocity component $u_{\textrm{ext}}^d$ with $d = \{1,2\}$:
\begin{subequations}
\begin{align}
	\grad{\levset_{\textrm{evo}}} \cdot \grad{u_{\textrm{ext}}^d} &= 0 \quad \quad \textrm{ in } \domain \in \gridNB, \label{eq:extensionProblemNonCut}\\
	u_{\textrm{ext}}^d &= u_{\interface}^d \quad \quad \textrm{ on } \interface.
\end{align}
\label{eq:extensionProblem}
\end{subequations}
In words, the extension velocity field $\velocExt$, resp. its components $u_{\textrm{ext}}^d$, is constructed by propagating the interface velocity $\velocI$ from $\interface$ perpendicular into the bulk domain $\bulk \in \gridNB$. 

The presented evolution algorithm is divided into two stages. In the first stage the extension velocity is computed monolithic with high-order accuracy on the cut-cells $\gridCC$ using the reformulation of the EVP \eqref{eq:extensionProblem} into an elliptic PDE. In this case the solution describes the steady-state limit of an anisotropic diffusion problem, see Utz and Kummer \cite{utz_high-order_2018}. In our implementation we follow the work of Utz\cite{utz_level_2018} and normalize the level-set gradient in \eqref{eq:extensionProblemNonCut} with $\grad{\phiDG}/\abs{\grad{\phiDG}}$ leading to more stable results. The discretization is done by the standard SIP method. The corresponding linear problem in $\velocExt$ is solved by a linear solver.

In the second stage the extension velocity on the near field $\gridNear$ is constructed geometrically in each cell separately using an adapted fast-marching\cite{adalsteinsson_fast_1995, adalsteinsson_fast_1999} procedure, see Algorithm \ref{alg:fastMarchExtVel}. For the low-order local solver in line \ref{algLine:localSolveExtVel} boundary values are provided on edges that are adjacent to cells in \textit{Accepted}. The solution is then constructed node-wise from the accepted edges into the cell $\cell_l$ and projected on $\mathbb{P}_k(\{K_l\})$. Note that Algorithm \ref{alg:fastMarchExtVel} is executed twice, once for negative near cells with $\phiDG < 0$ and once for $\phiDG > 0$.
\begin{algorithm}
		\caption{Fast-marching procedure for the EVP on near cells}
		\label{alg:fastMarchExtVel}
		\begin{algorithmic}[1]
			\Procedure{FastMarchingExtensionVelocity}{$\velocExt^{\textrm{cc}}$, $\gridCC$, $\gridNear$, $\phiDG$} \Comment{it is assumed that $\abs{\grad \phiDG} = 1$}
			\State $\velocExt^{\textrm{near}} \gets \velocExt^{\textrm{cc}}$; 
			\textit{Accepted} $\gets \gridCC$;
			\textit{Todo} $\gets \gridNear$
			\State \textit{Todo} $\gets \textrm{Sort}\left(\left| (1/\abs{\cell_l}) \int_{\cell_l} \phiDG \dif{V} \right| \leq \left| (1/\abs{\cell_{l+1}}) \int_{\cell_{l+1}} \phiDG \dif{V}\right| \right) \forall \cell_{l} \in$ \textit{Todo}, $\textrm{ for } 1 \leq l < L$ 
			\For{$l = 1, ..., L$}
			\State $\velocExt^{\textrm{near}+1} \gets \textrm{LocalSolve}\left(\cell_l, \textrm{Neighbor}(\cell_l) \in \right.$\textit{Accepted} $\left. \right)$ \Comment{solve Eq. \eqref{eq:extensionProblemNonCut} node-wise and project solution on $\cell_l$} \label{algLine:localSolveExtVel}
			\State \textit{Accepted} $\gets$ \textit{Accepted} $\cup \{K_l\}$; \textit{Todo} $\gets$ \textit{Todo} $\setminus \{K_l\}$
			\State $\velocExt^{\textrm{near}} \gets \velocExt^{\textrm{near}+1}$
			\EndFor
			\State \textbf{return} $\velocExt^{\textrm{near}}$
			\EndProcedure
		\end{algorithmic}
\end{algorithm}

The entire two-staged level-set evolution algorithm is presented in Algorithm \ref{alg:marchingLevSetEvo}. Since the signed-distance property of $\phiDG$ is only given on the narrow-band, a reinitialization is applied on new cells entering the narrow-band during the last advection. Therefore, a fast-marching procedure (see Algorithm \ref{alg:fastMarchReInit}) is executed in line \ref{algLine:ReInit}. After constructing the extension velocity field $\velocExt$ on the narrow-band $\gridNB$ (line \ref{algLine:EllipticExtVel} and \ref{algLine:FastMarchExtVel}) the new level-set field $\phiDG^{n+1}$ is computed according to the evolution equation \eqref{eq:levsetEvo} by 
\begin{equation}
\begin{aligned}
    \int_{\gridNB} \partial_t \phiDG \vert_{t^n} \vartheta \dif{V} - \int_{\gridNB} \left( \gradH{\phiDG^{n+1}} \cdot \velocExt \right) \vartheta \dif{V}\\ 
    + \oint_{\edge_{\textrm{int},\textrm{nb}}} \left( f_{\textrm{upwind}}\jump{\vartheta} - [\![\phiDG^{n+1} (\velocExt \cdot \normal_{\edge_{\textrm{int}, \textrm{nb}}}) \vartheta ]\!] \right) \dif{S} = 0, \quad \forall \vartheta \in \brknPspacek(\gridNB),
    \label{eq:levsetEvoDiscret}
\end{aligned}
\end{equation}
where $\edge_{\textrm{int}, \textrm{nb}}$ denotes the internal edges in $\gridNB$. The stabilization term $f_{\textrm{upwind}}$ describes a simple upwind-flux for the flux $\phiDG \velocExt$. The temporal term is discretized via the total variation diminishing Runge-Kutta scheme of third order according to Gottlieb and Shu\cite{gottlieb_total_1998}. At last, the advected level-set field $\phiDG^{n+1}$ is additionally stabilized by penalizing jumps at inner edges, i.e. $\int_{\edge_{\textrm{int}, \textrm{nb}}} (\pnlty_{\textrm{evo}}/\grdSz) \jump{\phiDG^{n+1}} \jump{\vartheta} \dif{S}$ with $\pnlty_{\textrm{evo}} = 10$, via an implicit Euler scheme with $\dt_{\textrm{pnlty}} = 0.001 \dt$. Thus, the deviation to the projected continuous level-set field $\levset$ used for the spatial discretization is reduced. 
\begin{algorithm}
	\caption{A two-staged marching level-set evolution algorithm on the narrow-band}
	\label{alg:marchingLevSetEvo}
	\begin{algorithmic}[1]
		\Procedure{MarchingLevelSetEvolution}{$\phiDG^{n}$, $\velocEvo$, $\dt$}
		\State $\grad{\phiDG^{n}} \gets \textrm{ComputeGradient}(\phiDG^{n}, \gridNear^{n})$
		\If{$\norm{\grad{\phiDG^{n}}} < 10^{-12}$ on $\gridNear^{n}$} \Comment{initialization of new cells in narrow-band}
		\State $\phiDG^{\ast} \gets \textrm{FastMarchReInit}(\phiDG^{n}, \gridCC^{n}, \gridNear^{n})$ \Comment{Algorithm \ref{alg:fastMarchReInit}} \label{algLine:ReInit}
		\EndIf
		\State $\grad{\phiDG^{\ast}} \gets \textrm{ComputeGradient}(\phiDG^{\ast}, \gridNear^{n})$
		\State $\velocExt^{\ast} \gets \textrm{EllipticExtensionVelocity}(\velocEvo, \gridCC^{n}, \phiDG^{\ast})$ \Comment{monolithic solve of elliptic formulation of EVP \eqref{eq:extensionProblem}} \label{algLine:EllipticExtVel}
		\State $\velocExt \gets \textrm{FastMarchingExtensionVelocity}(\velocExt^{\ast}, \gridCC^{n}, \gridNear^{n}, \phiDG^{\ast})$ \Comment{Algorithm \ref{alg:fastMarchExtVel}} \label{algLine:FastMarchExtVel}
		\State $\phiDG^{n+1} \gets \textrm{Advect}(\velocExt, \dt, \gridNB^{n})$ 
        \Comment{Eq. \eqref{eq:levsetEvoDiscret} Runge-Kutta scheme (TVD3)}
		\State $\phiDG^{n+1\ast} \gets \textrm{JumpPenalization}(\phiDG^{n+1})$ \Comment{via implicit Euler}
		\State $\phiDG^{n+1\ast} \gets \textrm{FullReInit}(\phiDG^{n+1\ast}, \gridNB^{n})$ 
		\Comment{optional for a predefined interval} \label{algLine:FastMarchReInit}
		\State \textbf{return} $\phiDG^{n+1\ast}$
		\EndProcedure
	\end{algorithmic}
\end{algorithm}

\paragraph{Level-set reinitialization}

As mentioned before the signed-distance property needs to be initialized for new cells entering the narrow-band (line \ref{algLine:ReInit} in Algorithm \ref{alg:marchingLevSetEvo}). Such a reinitialization problem is described by the following Eikonal equation
\begin{subequations}
\begin{align}
    \abs{\grad \phiDG} - 1 &= 0 \quad \textrm{ in } \domain,\\
    \phiDG &= 0 \quad \textrm{ on } \interface.
\end{align}
\label{eq:Eikonal}
\end{subequations}
Therefore, a fast-marching\cite{adalsteinsson_fast_1995, adalsteinsson_fast_1999} procedure is utilized as given in Algorithm \ref{alg:fastMarchReInit}. The local solver for one cell $K_j$ in line \ref{algLine:localSolveReInit} is based on the elliptic reinitialization method by Basting and Kuzmin\cite{basting_minimization-based_2013} and adapted for DG in Utz\cite{utz_interface-preserving_2017}. A low-order initial solution for $\phiDG$ is obtained node-wise by a direct geometric approach. Furthermore, 
note that at the end (line \ref{algLine:FastMarchReInit}) of the marching algorithm \ref{alg:marchingLevSetEvo} an optional full reinitialization of the level-set field $\phiDG$ on the narrow-band may be performed. In this case the elliptic reinitialization\cite{utz_interface-preserving_2017} is solved monolithic on the cut-cells and the fast-marching procedure \ref{alg:fastMarchReInit} is executed on the near field. 
\begin{algorithm}
		\caption{Fast-marching procedure for the reinitialization problem on near cells}
		\label{alg:fastMarchReInit}
		\begin{algorithmic}[1]
			\Procedure{FastMarchingReinitialization}{$\phiDG^{\textrm{cc}}$, $\gridCC$, $\gridNear$}
			\State \textit{Accepted} $\gets \gridCC$;
			\textit{Close} $\gets \textrm{Neighbour}(\gridCC) \in \gridNear$;
			\textit{Far} $\gets \cell_j \in \gridNear \setminus$ (\textit{Accepted} $\cup$ \textit{Close})
			\State $\phiDG^{\textrm{near}} \gets \phiDG^{\textrm{cc}}$
			\While{\textit{Close} $\neq \emptyset$}
			\State $\phiDG^{\textrm{near}+1} \gets \textrm{LocalSolve}(\phiDG^{\textrm{near}}) \forall \cell_j \in$ \textit{Close} \Comment{Solve Eq. \eqref{eq:Eikonal}} \label{algLine:localSolveReInit}
			\State \textit{Accepted} $\gets \cell_j: \min \left( \left| (1/\abs{\cell_j}) \int_{\cell_j} \phiDG^{n+1} \dif{V} \right| \right)$ \Comment{Remove $\cell_j$ from \textit{Close}}
			\State \textit{Close} $\gets \textrm{Neighbor}\left( \cell_j: \min \left(\left| (1/\abs{\cell_j}) \int_{\cell_j} \phiDG^{n+1} \dif{V} \right| \right) \right) \in \textit{Far}$ \Comment{Remove neighbors from \textit{Far}}
			\State $\phiDG^{\textrm{near}} \gets \phiDG^{\textrm{near}+1}$
			\EndWhile
			\State \textbf{return} $\phiDG^{\textrm{near}}$
			\EndProcedure
		\end{algorithmic}
\end{algorithm}


\section{The XDG method for the instationary two-phase NSE}
\label{sec:XDGmethod}

In this section the XDG discretization of the considered transient two-phase flow problem is provided (Section \ref{sec:XDGdiscretization}) and the concept of the cell agglomeration for resolving stability issues and handling topology changes during interface motion is described (Section \ref{sec:CellAgglom}).

The formulation of an XDG method requires the numerical integration of cut-cell integrals such as
\begin{equation}
    \int_{\cell_{\spc}^{\cut}} f \dif{V}, \quad \int_{\partial\cell_{\spc}^{\cut}} f \dif{S} \quad \textrm{and} \quad \int_{\interface} f \dif{S},
    \label{eq:cutcellIntegrals}
\end{equation}
which represent integrals over cut-cell volumes and surfaces and integrals along the interface $\interface$. An essential prerequisite of the XDG method is a high-order accurate integration technique allowing for a stable and robust numerical integration on implicitly defined domains.
In contrast to Kummer\cite{kummer_extended_2016} we use in this work the quadrature method proposed by Saye \cite{saye_high-order_2015}. It is generally faster compared to the Hierarchical Moment Fitting, but restricted to hyperrectangles.

\subsection{XDG discretization of the incompressible two-phase NSE}
\label{sec:XDGdiscretization}

Before introducing the discretization for the transient two-phase flow problem, we define the function space of the ansatz and test functions for the velocity field $\veloc \in \mathbb{R}^2$ and the pressure field $\press \in \mathbb{R}$ by
\begin{equation}
	\spaceVkX(t) := \prod_{k_{\gamma}}^{\vectr{k}} \mathbb{P}_{k_\gamma}^{\cut}(\grid, t),
	\label{eq:functionSpace}
\end{equation}
where  $\vectr{k} = \{ k, k, k^{\prime} \} = k_{\gamma}$ with $\gamma = 1,2,3$ describes the degree vector. Here, the velocity components are discretized by an XDG space of order $\Pdeg$ and the pressure field of order $\Pdeg^{\prime} = \Pdeg - 1$. Thus, we comply with the Lady\u{z}enskaja-Babu\u{s}ka-Brezzi condition \cite{babuska_finite_1973, brezzi_existence_1974} for the Stokes system with $\indA{\dens} = \indB{\dens}$ and $\indA{\visc} = \indB{\visc}$. 

Following the spatial discretization in Kummer \cite{kummer_extended_2016}, we propose the following discretization for the transient two-phase incompressible Navier-Stokes equations \eqref{eq:NavierStokes} with the interface jump conditions \eqref{eq:jump_NavierStokes} and boundary conditions \eqref{eq:boundary_conditions}:
Find $(\veloc^{n+1}, \press^{n+1}) \in \spaceVkX(t^{n+1})$, such that $\forall (\testV^{n+1}, \testP^{n+1}) \in \spaceVkX(t^{n+1})$
\begin{equation}
\begin{aligned}
	m^{n+1}(\veloc^{n+1}, \testV^{n+1}) + d_0 \dt \left( c(\veloc^{n+1}, \veloc^{n+1}, \testV^{n+1}) + b(\press^{n+1}, \testV^{n+1}) \right.\\ \left.
	- a(\veloc^{n+1}, \testV^{n+1}) - b(\testP^{n+1}, \veloc^{n+1}) \right) = d_0 \dt \ s(\testV^{n+1},\testP^{n+1})\\ 
	+ d_1 m^{n}(\veloc^{n}, \testV^{n}) + d_2 m^{n-1}(\veloc^{n-1}, \testV^{n-1}) + d_3 m^{n-2}(\veloc^{n-2}, \testV^{n-2}).
\label{eq:variationalForm_NSE} 
\end{aligned}
\end{equation}
The bilinear forms denoted by $m^n(-,-)$ represent the mass matrix components for the corresponding time step $t^{n}$ with
\begin{equation}
	m^n(\veloc^{n}, \testV^{n}) = \int_{\domain \setminus \interface(t^n)} \dens \veloc^{n} \cdot \testV^{n} \dif{V}.
\label{eq:massTerm}
\end{equation}
A third order backward differentiation formula (BDF) is employed, where the coefficients for each time step are given by $d_0 = 6/11$, $d_1 = 18/11$, $d_2 = -9/11$, $d_3 = 2/11$. One should note that in the context of the moving interface approach by Kummer et al. \cite{kummer_time_2018} a spatial approximation of degree $\Pdeg$ theoretically requires at least a time integration scheme of order $2\Pdeg$ for linear equations. BDF schemes are not $A$-stable above order 2 including more eigenvalues $z$ with $\Re(z) < 0$ for higher orders. However, the unstable eigenvalues for the third order scheme are comparable small, i.e. near the imaginary axis. Since we are dealing with comparably small time steps
due to the capillary time step restriction (see Eq. \eqref{eq:capillaryTimestep}), we choose the use of the third order BDF scheme to comply with the requirement for the moving interface approach. 

In the remainder of this section the integration domains of the remaining terms are set to the time-level $t^{n+1}$, i.e. the interface is fixed at $\interface = \interface(t^{n+1})$. The trilinear form $c(\veloc^{\ast}, \veloc, \testV)$ describes the discretization of the convective term, where a local Lax-Friedrichs flux is employed
\begin{equation}
	c(\veloc^{\ast}, \veloc, \testV) = - \int_{\domain} \dens (\veloc \otimes \veloc^{\ast}) : \gradH{\veloc} \dif{V}
	- \oint_{\edgeInt \cup \edgeN \setminus \interface} \left( \aver{\veloc \otimes \veloc^{\ast}} \normalGam + \frac{\lambda}{2} \jump{\veloc} \right) \cdot \jump{\dens \testV} \dif{S}.
	\label{eq:convectiveTerms}
\end{equation}
For the choice of the stabilization parameter $\lambda$ we refer to Kummer \cite{kummer_extended_2016}. Note that the interface integral is excluded from the surface discretization. This contribution is canceled out by 
\begin{equation}
    \int_{\interface} (\veloc^{\ast} \cdot \normalI) \veloc \dif{S},
\end{equation}
which is introduced within the moving interface approach\cite{kummer_time_2018}. The discretized problem \eqref{eq:variationalForm_NSE} needs to be solved iteratively due to the non-linearity of the convective term. Thus, the whole system is linearized in each iteration with $\veloc^{\ast}$. The iteration process in combination with the interface evolution is presented in Section \ref{sec:Solver}. 

The bilinear form $b(-,-)$ represents both, the discretization of the pressure gradient and the continuity term, i.e. velocity divergence,
\begin{equation}
	b(\press, \testV) = - \int_{\domain} \press \divH{\testV} \dif{V} - \oint_{\edgeInt \cup \edgeD} \jump{\testV} \cdot \normalGam \aver{\press} \dif{S}.
	\label{eq:pressureContiTerms}
\end{equation}
An extension to the standard symmetric interior penalty (SIP) method \cite{arnold_interior_1982} is used for the discretization of the viscous terms. In this case the transposed term $\gradT{\veloc}$ within the divergence operator is included and since $\gradHT{\veloc} : \gradH{\testV} = \gradHT{\testV} : \gradH{\veloc}$ the bilinear form $a(\veloc, \testV)$ is still symmetric in $\veloc$ and $\testV$ with
\begin{equation}
\begin{aligned}
	a(\veloc, \testV) = &-\int_{\domain} \visc \left( \gradH{\veloc} : \gradH{\testV} + \gradHT{\veloc} : \gradH{\testV} \right) \dif{V} \\
	&+ \oint_{\edgeInt \cup \edgeD} \left( \aver{ \visc \left( \gradH{\veloc} + \gradHT{\veloc} \right) } \normalGam \right) \cdot \jump{\testV} \dif{S}\\
    + \oint_{\edgeInt \cup \edgeD} &\left( \aver{ \visc \left( \gradH{\testV} + \gradHT{\testV} \right) } \normalGam \right) \cdot \jump{\veloc} \dif{S}
	- \oint_{\edgeInt \cup \edgeD} \pnlty \jump{\veloc} \cdot \jump{\testV} \dif{S}.
\label{eq:viscousTerms}
\end{aligned}
\end{equation}
The corresponding penalty parameter $\pnlty$ is chosen according to 
\begin{align}
	\pnlty := \begin{cases}
	\max{ \{ \inn{\visc}, \out{\visc} \} } \max{ \{ \inn{\tilde{\pnlty}}, \out{\tilde{\pnlty}} \} } \quad &\text{ on } \edgeInt,\\
	\inn{\visc} \inn{\tilde{\pnlty}} \quad &\text{ on } \boundary,
	\end{cases}
\label{eq:pnltyParam}
\end{align}
where the local penalty parameter $\tilde{\pnlty}$ is computed by 
\begin{equation}
	\tilde{\pnlty} = \pnlty_0 \Pdeg^2 \frac{\abs{\partial{\cell_{\spc}^{\cut}}}}{\abs{\cell_{\spc}^{\cut}}}.
\label{eq:pnltyFactor}
\end{equation}
The safety factor $\pnlty_0$ is set to $\pnlty_0 = 4$, see Kummer \cite{kummer_extended_2016} for details. The geometric factor $\abs{\partial{\cell_{\spc}^{\cut}}} / \abs{\cell_{\spc}^{\cut}}$ is directly provided by the construction of the cut-cell quadrature rules, see Eq. \eqref{eq:cutcellIntegrals}. 

Finally, we specify the term $s(\testV, \testP)$ on the right-hand side of the variational formulation \eqref{eq:variationalForm_NSE}, which summarizes the Dirichlet boundary conditions and force terms
\begin{equation}
\begin{aligned}
	s(\testV, \testP) = &- \oint_{\edgeD} \dens \left( (\diri{\veloc} \otimes \diri{\veloc}) \normalGam + \frac{\lambda}{2} \diri{\veloc} \right) \cdot \testV \dif{S} -\oint_{\edgeD} \diri{\veloc} \cdot \left( \visc \left( \gradH{\testV} + \gradHT{\testV} \right) \normalGam - \pnlty \testV \right) \dif{S}\\
	&+ \oint_{\edgeD} \testP \diri{\veloc} \cdot \normalGam \dif{S} +\int_{\domain} \dens \vectr{g} \cdot \testV \dif{V} 
	- \oint_{\interface} \surfT \projI : \gradI{\testV} \dif{S} 
	+ \int_{\interface \cap \edge} \surfT \aver{\tangent} \cdot \jump{\testV} \dif{l}.
\label{eq:rhsTerms}
\end{aligned}
\end{equation}
The first two terms in the upper line describe the Dirichlet boundary conditions for the convective part \eqref{eq:convectiveTerms} and the viscous parts \eqref{eq:viscousTerms}. The first term in the bottom line corresponds to the Dirichlet boundary condition for the continuity equation \eqref{eq:pressureContiTerms} and the second term denotes the discretization of the volume term $\dens \vectr{g}$. 

For the numerical treatment of the surface tension force (remaining two terms in the second line) we do not evaluate the curvature $\curv$, see \eqref{eq:jump_momentumMaterial}, but choose the Laplace-Beltrami formulation, see e.g. Gross and Reusken \cite{gross_finite_2007}, Cheng and Fries \cite{cheng_xfem_2012} and Sauerland and Fries \cite{sauerland_stable_2013}. In this case the surface tension force is rewritten to its interface divergence form with $\surfT \curv \normalI = \divI{\left( \surfT \projI \right)}$, where $\projI := \tensr{I} - \normalI \otimes \normalI$ describes the projection tensor onto the interface $\interface$ and the interface gradient operator $\gradI{}$ is defined as $\gradI{} := \projI \nabla$. The vector $\tangent$ denotes the tangential to $\interface$ and the corresponding integral denotes point-measures in $\mathbb{R}^2$. Those boundary terms ensure force conservation of the surface tension along closed interfaces. Note that we do not employ a semi-implicit discretization as introduced in Dziuk\cite{dziuk_algorithm_1990} or a regularization within the continuum surface force framework\cite{brackbill_continuum_1992, bansch_finite_2001, hysing_new_2006}.

\subsection{Cell agglomeration}
\label{sec:CellAgglom}

The discretization by the XDG method may introduce arbitrarily small cut-cells $\cell_{j,\spc}^{\cut}$, where the volume fraction $\abs{\cell_{j,\spc}^{\cut}} / \abs{\cell_{j}}$ is small compared to the background cell $\cell_{j}$. In context of the SIP discretization for the viscous terms \eqref{eq:viscousTerms} such cut-cells give rise to undesirably high condition numbers and consequently leads to stability issues. Therefore, cell agglomeration is employed to remove small cut-cells from the discretized system\cite{kummer_extended_2016}. Furthermore, the cell agglomeration is used in the context of the moving interface time discretization \cite{kummer_time_2018}. In this case, appearing and vanishing cut-cells are handled during the interface motion.

\begin{figure}
	\centering
	\def\svgwidth{400pt}
	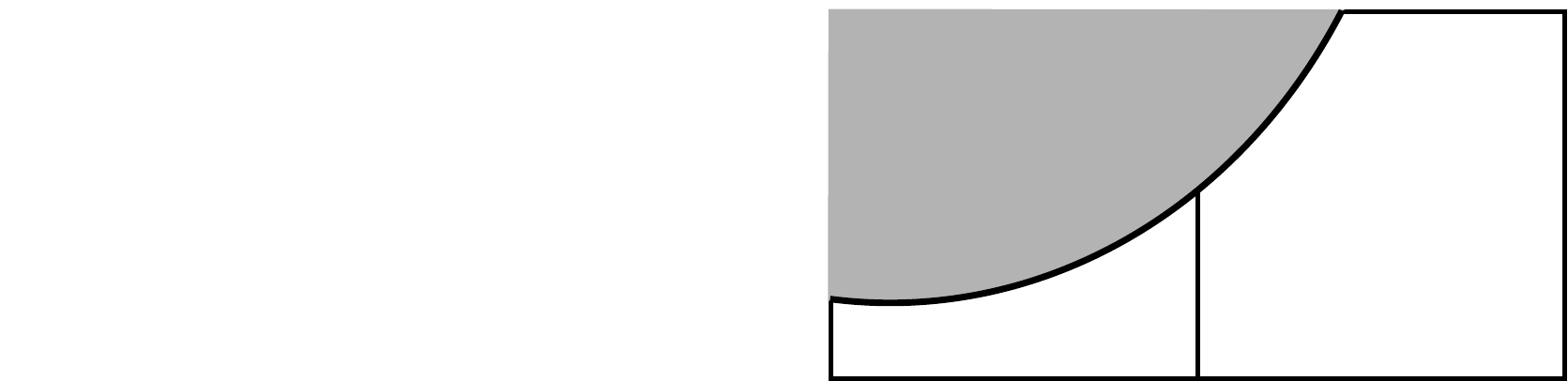
	\caption{Left: Definition of a cut-cell mesh $\grid^{\cut}$ as an undirected graph $G(\grid^{\cut}) := (\grid^{\cut}, \logEdg(\grid^{\cut}))$ by logical edges $\logEdg(\grid^{\cut}) = \{ \logEdg_1(\{ \cell_{j,\domainA}^{\cut}, \cell_{j,\domainB}^{\cut}\}), \logEdg_2(\{ \cell_{j,\domainA}^{\cut}, \cell_{g,\domainA}^{\cut}\}), \logEdg_3(\{ \cell_{j,\domainB}^{\cut}, \cell_{g,\domainB}^{\cut}\}), \logEdg_4(\{ \cell_{g,\domainA}^{\cut}, \cell_{g,\domainB}^{\cut}\}) \}$. Right: Aggregation mesh $\Agg(\grid^{\cut}, \AggMap)$ of the aggregation map $\AggMap = \{ \logEdg_2 \}$.}
	\label{fig:cutCellAgg}
\end{figure}

The cell agglomeration and resulting meshes are described in terms of graph theory according to Kummer et al. \cite{kummer_bosss_2020}. Basic definitions for the agglomeration are given in the following. An example is provided is in Figure \ref{fig:cutCellAgg}.
\begin{definition}[logical edge and cell agglomeration] For a cut-cell mesh $\grid^{\cut}$ we define:
\begin{itemize}
    \item{a logical edge $\logEdg(\{ \cell_{j,\spc}^{\cut}, \cell_{g,\spc}^{\cut}\})$ for two neighboring cut-cells $\cell_{j,\spc}^{\cut}$ and $\cell_{g,\spc}^{\cut}$.} 
    \item{an undirected graph $G(\grid^{\cut}) := \left( \grid^{\cut}, \logEdg(\grid^{\cut}) \right)$, where $\logEdg(\grid^{\cut})$ denotes all logical edges in the cut-cell mesh $\grid^{\cut}$.}
    \item{a cluster of neighboring cells $a = \{ \cell_{l,\spc}^{\cut}, ..., \cell_{L,\spc}^{\cut} \}$ and a corresponding aggregated cell $\cell_{a,\spc}^{\cut} :=  \cup_{\cell_{l,\spc}^{\cut} \in a}^{\partial} \cell_{l,\spc}^{\cut}$. The modified union is defined as $ X \cup^{\partial} Y :=  (\overline{X} \cup \overline{Y}) \setminus \partial(\overline{X} \cup \overline{Y})$ in order to ensure simple connectivity of $\cell_{a,\spc}^{\cut}$. All elements in $a$ are related to the same species.} 
    \item{an aggregation mesh $\Agg(\grid^{\cut}, \AggMap)$ including all aggregated cells given by an aggregation map $\AggMap \subset \logEdg(\grid^{\cut})$ and including non-aggregated cells.}
\end{itemize}
\end{definition}

\begin{definition}[agglomerated XDG space] For some aggregation map $\AggMap \subset \logEdg(\grid^{\cut})$ and corresponding aggregation mesh $\Agg(\grid^{\cut}, \AggMap)$ the agglomerated broken cut-polynomial space of total degree $\Pdeg$ is defined as:
\begin{equation}
	\brknPspacek^{\cut, \AggMap}:=  \brknPspacek(\Agg(\grid^{\cut}, \AggMap)).
\end{equation}
Note that $\brknPspacek(\Agg(\grid^{\cut}, \AggMap))$ is a sub-space of the original XDG space $\brknPspacek(\grid^{\cut})$. A corresponding basis to the aggregation mesh $\Agg(\grid^{\cut}, \AggMap)$ is denoted by ${\vectr{\basis}}^{\cut, \AggMap}$.
\end{definition}

\paragraph{Agglomeration for removing small cut-cells} 

In order to remove small cut-cells form the numerical mesh $\grid^{\cut}$, we introduce an agglomeration map $\AggMap_{\AggTh} \subset \logEdg(\grid^{\cut})$ with edges $\logEdg(\{ \cell_{j,\spc}^{\cut}, \cell_{g,\spc}^{\cut}\})$ matching the following two conditions:
\begin{itemize}
	\item{The volume fraction of $\cell_{j,\spc}^{\cut}$ is smaller than the agglomeration threshold $\AggTh$, i.e $0 < \abs{\cell_{j,\spc}^{\cut}} / \abs{\cell_{j}} < \AggTh$.}
	\item{The cell $\cell_{g,\spc}^{\cut}$ marks the neighbor with the largest volume fraction in the same species.}
\end{itemize}
In this work the agglomeration threshold $\AggTh$ is set to $\AggTh = 0.1$ for all presented simulations. A condition number study with varying agglomeration thresholds is presented in Kummer\cite{kummer_extended_2016}.

\paragraph{Agglomeration for topology changes during interface motion}

\begin{figure}
	\centering
	\def\svgwidth{300pt}
	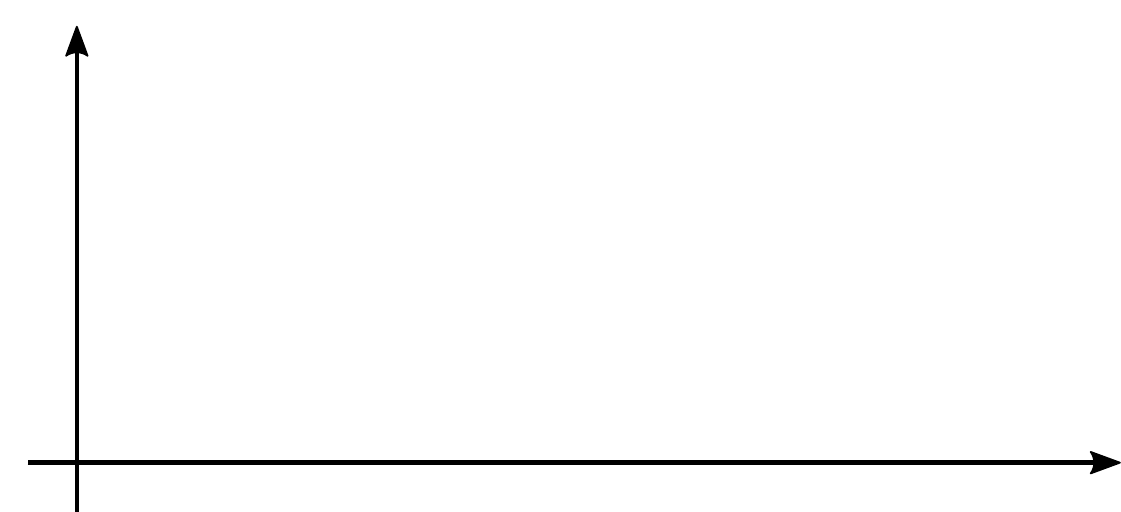
	\caption{The interface $\interface$ crosses a background cell edge during on time step resulting in a vanishing cut-cell $\cell_{j,\domainB}^{\cut}$, i.e. $\abs{\cell_{j,\domainB}^{\cut}(t^0)} > 0$ and $\abs{\cell_{j,\domainB}^{\cut}(t^1)} = 0$, and an appearing cut-cell $\cell_{g,\domainA}^{\cut}$, i.e. $\abs{\cell_{g,\domainA}^{\cut}(t^0)} = 0$ and $\abs{\cell_{g,\domainA}^{\cut}(t^1)} > 0$. Thus, the agglomeration maps $\AggMap(t^0) = \{ \logEdg_\textrm{van} \}$ and $\AggMap(t^1) = \{ \logEdg_\textrm{new} \}$ ensure equal topology on both time steps $t^0$ and $t^1$.}
	\label{fig:cutCellAggTime}
\end{figure}

Ensuring the same topology throughout all mass matrices \eqref{eq:massTerm} at the different time steps, appearing and vanishing cut-cells due to the interface motion are handled via cell agglomeration, see Figure \ref{fig:cutCellAggTime}. Therefore, the agglomeration map $\AggMap_{\AggTh}$ additionally exhibits logical edges $\logEdg(\{ \cell_{j,\spc}^{\cut}, \cell_{g,\spc}^{\cut}\})$ matching the conditions:
\begin{itemize}
	\item{The cut-cell $\cell_{j,\spc}^{\cut}$ appears, i.e. $\abs{\cell_{j,\spc}^{\cut}(t^{n})} = 0$ and $\abs{\cell_{j,\spc}^{\cut}(t^{n+1})} > 0$, or vanishes, i.e. $\abs{\cell_{j,\spc}^{\cut}(t^{n})} > 0$ and $\abs{\cell_{j,\spc}^{\cut}(t^{n+1})} = 0$.}
	\item{The cell $\cell_{g,\spc}^{\cut}$ marks the neighbor with the largest volume fraction in the same species.}
\end{itemize}

%
\section{Solver structure -  Coupling flow solver and interface evolution}
\label{sec:Solver}

In this section the coupling of the flow solver for the discretizied system \eqref{eq:variationalForm_NSE} and the interface evolution is presented. Additionally, the adaptive mesh refinement (AMR) procedure is briefly discussed at the end.

In Algorithm \ref{alg:RunSolverOneStep} the procedure for a single solver run is given. There are two options for the coupling of the interface evolution: explicit and implicit. First, we describe the explicit coupling, where the interface is only updated once at the beginning (Line \ref{algLine:UpdateLS}). An initial solution on the new domain $\Agg(\grid^{\cut}, \AggMap_{\AggTh}^{n+1})$ is given by extrapolation of the old solutions $\veloc^{n}$ and $\press^{n}$. The extrapolation operator is defined as follows
\begin{equation}
\begin{aligned}
	\textrm{Ex}_{\left( t^n, t^{n+1} \right)}: \brknPspacek\left(\Agg\left(\grid^{\cut},\AggMap^n\right) \right) &\rightarrow \brknPspacek\left(\Agg\left(\grid^{\cut},\AggMap^{n+1}\right) \right)\\
	 \vectr{\tilde{f}}^n \cdot \vectr{\basis}^{\cut, \AggMap^n}\left( t^n \right) &\mapsto \vectr{\tilde{f}}^n \cdot \vectr{\basis}^{\cut, \AggMap^{n+1}}\left( t^{n+1} \right).
\end{aligned}
\label{eq:extrapolation}
\end{equation}
Note that the coefficients $\vectr{\tilde{f}}^n$ do not change. This initial solution further provides a first linearization $\veloc^{\ast}$ (Line \ref{algLine:UpdateLin}) for the convective terms \eqref{eq:convectiveTerms}. After that, the procedure enters the iterative solution process via Picard iterations until the prescribed convergence criterion $\ccNSE$ is satisfied or a maximum number of iterations $i_\textrm{max}$ is exceeded (Line \ref{algLine:Iteration}).
\begin{algorithm}
	\caption{One time step of the XNSE-solver}
	\label{alg:RunSolverOneStep}
	\begin{algorithmic}[1]
		\Procedure{RunSolverOneStep}{$t^n$,$\dt$} \Comment{solve $\levset$, $\veloc$ and $\press$ for $t^{n+1} = t^n + \dt$}
		\State $\levset^{n}, \veloc^{n}, \press^{n} \gets \levset^{n+1}, \veloc^{n+1}, \press^{n+1}$ \Comment{set previous time step to outdated $t^{n+1}$} 
		\State $\levset^{n+1,0} \gets $ \textproc{UpdateLevelSet}($\veloc^{n}$, $\dt$) \Comment{Algorithm \ref{alg:UpdateLevelSet}} \label{algLine:UpdateLS}
		\State $\Agg\left(\grid^{\cut}, \AggMap_{\AggTh}^{0}\right) \gets \grid^{\cut}$ \Comment{update agglomeration}
		\State $\massM^{n+1,0} \gets \textrm{ComputeMassMatrix}()$ 
		\State $\veloc^{n+1,0}, \press^{n+1,0} \gets \textrm{Extrapolation}(\veloc^{n}, \press^{n}, \levset^{n+1,0}, \levset^{n})$ \Comment{Initial solution on $\levset^{n+1,0}$ by Eq. \eqref{eq:extrapolation}} \label{algLine:Extrapolation}
		\State $\veloc^{\ast} \gets \textrm{UpdateLinearization}(\veloc^{n+1,0})$ \Comment{see convective terms Eq. \eqref{eq:convectiveTerms}} \label{algLine:UpdateLin}
		\State $\res^0 \gets \textrm{ComputeResidual}(\veloc^{n+1,0}, \press^{n+1,0})$
		\While{(!($\res^i < \ccNSE$) $\&\&$ $i < i_\textrm{max}$)} \Comment{if implicit: !($\res^i < \ccNSE$ $\&\&$ $\res_\textrm{LS}^i < \ccLS$)} \label{algLine:Iteration}
		\If{(implicit interface evolution $\&\&$ $\res^i < \ccNSE$)} \label{algLine:Implicit}
		\State $\levset^{n+1,i} \gets $ \textproc{UpdateLevelSet}($\veloc^{n+1,i}$, $\dt$)
		\State $\Agg\left(\grid^{\cut}, \AggMap_{\AggTh}^{i}\right) \gets \grid^{\cut}$
		\State $\massM^{n+1,i} \gets \textrm{ComputeMassMatrix}()$
		\EndIf
		\State $\OpM^i, \RHS^i \gets \textrm{ComputeOperatorMatrixAndRHS}(\veloc^{n+1,i}, \press^{n+1,i})$ \Comment{for Eq. \eqref{eq:matrixForm_NSE}} \label{algLine:ComputeOpRHS}
		\State $\massM_{\AggMap}^{n+1,i}, \OpM_{\AggMap}^{i}, \RHS_{\AggMap}^{i} \gets \textrm{PerformAgglomeration}\left(\Agg\left(\grid^{\cut}, \AggMap_{\AggTh}^{i}\right)\right)$ \Comment{update old mass matrices} \label{algLine:PerformAgg}
		\State $\OpM_{\AggMap}^{i, \ast} \gets$ \textrm{Preconditioning}$\left( \OpM_{\AggMap}^{i} \right)$ \label{algLine:precondition} 
		\State $\veloc_{\AggMap}^{n+1,i}, \press_{\AggMap}^{n+1,i} \gets \textrm{SolveLinearSystem}()$ \Comment{solve Eq. \eqref{eq:matrixForm_NSE} by sparse direct solver}
		\State $\veloc^{n+1,i}, \press^{n+1,i} \gets \textrm{RevertAgglomeration}()$ \label{algLine:RevertAgg}
		\State $\veloc^{\ast} \gets \textrm{UpdateLinearization}(\veloc^{n+1,i})$ \Comment{Eq. \eqref{eq:convectiveTerms}} 
		\State $\res^i \gets \textrm{ComputeResidual}(\veloc^{n+1,i}, \press^{n+1,i})$
		\State $i \gets i+1$
		\EndWhile
		\State \textrm{PostProcessing}($t^n$,$\dt$)
		\EndProcedure
	\end{algorithmic}
\end{algorithm}

Within each iteration the linear system of the discretizied Navier-Stokes equations \eqref{eq:variationalForm_NSE} is updated (Line \ref{algLine:ComputeOpRHS}). The corresponding matrix formulation is written as
\begin{equation}
    \left( \massM^{n+1} + d_0 \dt \ \OpM \right) \vectr{s}^{n+1} = d_0 \dt \ \RHS + d_1 \massM^{n} \vectr{s}^{n} + d_2 \massM^{n-1} \vectr{s}^{n-1} + d_3 \massM^{n-2} \vectr{s}^{n-2},
    \label{eq:matrixForm_NSE}
\end{equation}
with $\vectr{s}^{n+1} = (\tilde{\veloc},\tilde{\vectr{\press}})$ denoting the sought-after solution vector including the coefficients of the velocity fields $\tilde{\veloc}$ and the pressure field $\tilde{\vectr{\press}}$. The corresponding operator matrix $\OpM$ has the structure of a saddle point problem given by
\begin{equation}
	\OpM = \left[ \begin{array}{cc}
	\OpM_{c,a} & \OpM_b^\textrm{T} \\
	\OpM_b & 0
	\end{array} \right],
\end{equation}
where $\OpM_{c,a}$ denotes the convective trilinear form $c(-,-,-)$ in \eqref{eq:convectiveTerms} and viscous bilinear form $a(-,-)$ in \eqref{eq:viscousTerms}, and $\OpM_b$ the bilinear forms for the pressure gradient and velocity divergence $b(-,-)$ in \eqref{eq:pressureContiTerms}. After performing the cell agglomeration (Line \ref{algLine:PerformAgg}) the resulting system is preconditioned further reducing the condition number (Line \ref{algLine:precondition}). Therefore, a block Jacobi preconditioning is used resulting in diagonal matrices containing only $0$, $-1$ and $+1$ entries for the symmetric part of the block diagonal in $\OpM_{c,a}$, see Kummer\cite{kummer_extended_2016} for details. The preconditioned linear system is solved using the direct solver MUMPS\cite{amestoy_fully_2001, amestoy_hybrid_2006}.\\

The procedure described above corresponds to an explicit coupling, where the interface position is updated once and is determined by the velocity field of the previous time step, i.e. $\velocEvo = \veloc^{n}$. Employing an implicit coupling, i.e. $\velocEvo = \veloc^{n+1}$, the level-set update is additional performed within the iteration process (Line \ref{algLine:Implicit}). Note that in this case we do not update the level-set in every iteration, but when the Navier-Stokes solution is converged. The coupled iteration process converges when both the Navier-Stokes solution with $\ccNSE$ and the level-set field with $\ccLS$ converge.

\paragraph{Adaptive mesh refinement}

The XNSE-solver is extended to allow adaptive mesh refinement (AMR) during the simulation. Therefore, a background cell is divided into four equal sized sub-cells. 
For higher refinement levels we ensure that neighboring cell also refine in order to exhibit always a 2:1 cell ratio on every edge. Thus, we counteract undesired locking-effects on the refined cells by much larger cells. Furthermore, we do not coarsen on cut-cells.
The adaption of the mesh and the corresponding approximated DG-fields is performed before \textproc{RunSolverOneStep} and requires a soft-restart of the current simulation. The indication of the current refinement level of each cell is predefined. In general we employ a constant refinement on the narrow band, but additional refinement may be used depending on the flow characteristics. 

%
\section{Numerical results}
\label{sec:results}

The following simulations are all computed with an explicit coupling of the level-set evolution and no additional reinitialization is performed. Since approximating the interface, resp. the level-set field, by a polynomial basis of degree $\Pdeg$ we set the minimal resolved capillary wave length to $\frac{\grdSz}{\Pdeg+1}$. The resulting capillary time step according to Denner and van Wachem\cite{denner_numerical_2015} is given by 
\begin{equation}
	\dt_{\surfT} = \sqrt{\frac{({\dens}_{\domainA} + {\dens}_{\domainB}) \left({ \frac{\grdSz}{\Pdeg+1}}\right)^3}{2\pi\surfT}}. 
	\label{eq:capillaryTimestep}
\end{equation}
The numerical time step size $\dt$ always satisfies $\dt < \dt_{\surfT}$, if not stated otherwise. All simulations in this section are performed with the XNSE-Solver within the open-source solver framework BoSSS\cite{kummer_bosss_2020} developed at the department of fluid dynamics at the TU Darmstadt. The source code is available under the Apache License at \url{http://github.com/FDYdarmstadt/BoSSS}. The presented results are available online at \url{https://doi.org/10.25534/tudatalib-327}

In this results section we often derive scalar measures over the simulation time from the computed numerical DG fields. In order to quantify the error against a reference solution we define the $l_2$-error according to Hysing et al.\cite{hysing_quantitative_2009} by 
\begin{align}
	l_2 \textrm{-error}: \ltwoErr{e} &= \left(\frac{\sum_{n=1}^{\textrm{NTS}} \abs{q_n - q_{n,\textrm{ref}}}^2 }{\sum_{n=1}^{\textrm{NTS}} \abs{q_{n,\textrm{ref}}}^2 } \right)^\frac{1}{2}, \label{eq:ScalarErrorNorms2}
\end{align}
where $q_N$ describes the scalar measure at time step $t^n$, i.e. $q_n = q(t^n)$ with $n = 1, ..., \textrm{NTS}$ and NTS denoting the total number of time steps. $q_{n,\textrm{ref}}$ describes a suitable reference solution at $t^n$, which is either an analytical solution or the solution on the finest grid of a convergence study. If there are less time steps than provided by the reference solution, we use linear interpolation to generate missing data. We determine the experimental order of convergence (EOC) by the slope of a linear regression in the corresponding log-log plot. The regression coefficients are estimated by an ordinary least-squares method. Further, we define the rate of convergence (ROC) for a single refinement level $l$ according to Hysing et al.\cite{hysing_quantitative_2009} by 
\begin{equation}
	\textrm{ROC} = \frac{\log_{10}(\norm{e^{l-1}}/\norm{e^l})}{\log_{10}(\grdSz^{l-1}/\grdSz^l)}.
	\label{eq:ConvergenceRate}
\end{equation} 

In the following we present the results of capillary waves and comparing to the analytical solution for the amplitude height (Section \ref{sec:CW}). Simulations regarding a droplet in equilibrium and non-equilibrium state (Section \ref{sec:droplet}). At last we compare our solver against the benchmark groups of the rising bubble test case (Section \ref{sec:RB}).  

\subsection{Capillary wave}
\label{sec:CW}

\begin{figure}
	\centering
	\def\svgwidth{140pt}
	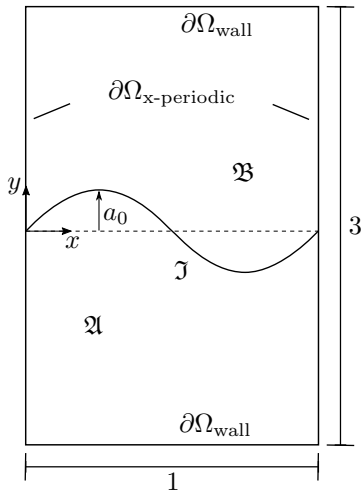
	\caption{Numerical setup for the capillary wave test case with an sinusoidal initial perturbation of wavelength $\lambda = 1$ and amplitude height $a_0 \ll \lambda$.}
	\label{fig:CWsetup}
\end{figure}

We consider the damped oscillations of a capillary wave initiated by a sinusoidal disturbance with a wavelength of $\lambda = 1$ and an initial amplitude height of $a_0 = 0.01 \ll \lambda$, see Figure \ref{fig:CWsetup}. The physical properties in both phases are set equal. For such a setting, i.e. small-amplitude limit and same kinematic viscosities $\nu = \visc / \dens$, Prosperetti \cite{prosperetti_motion_1981} presents an analytical solution of the amplitude height $a(t)$ for the corresponding initial value problem. 
Comparing to the analytical solution we follow the work of Popinet \cite{popinet_accurate_2009} and set the computational domain to $\domain = [0,\lambda] \times [-3\lambda/2, 3\lambda/2]$. The lower and upper boundaries are imposed by a wall boundary condition and the side boundaries are periodic in $x$-direction. In order to verify against a representative range of physical regimes, from an overdamped oscillation (low La) to a highly oscillatory behavior (high La), we investigate a study of the Laplace number defined as $\laplace = \surfT \rho \lambda / \visc^2$. The study consists of four values with $\laplace = \{ 3, 1.2 \cdot 10^{2}, 3 \cdot 10^{3}, 3 \cdot 10^{5}\}$ and the corresponding physical properties are given in Table \ref{tab:CWparamLaStudy}. 
\begin{center}
\begin{table}
\caption{Physical parameters of the Laplace number study for capillary waves}
\centering
	\begin{tabular}{c|cccc}
		$\laplace$ & $3$ & $1.2 \cdot 10^{2}$ & $3 \cdot 10^{3}$ & $3 \cdot 10^{5}$\\
		\hline 
		$\dens$ & $1\cdot 10^{-3}$ & $1\cdot 10^{-3}$ & $1\cdot 10^{-3}$ & $1\cdot 10^{-3}$ \\
		$\visc$ & $1\cdot 10^{-3}$ & $5\cdot 10^{-4}$ & $1\cdot 10^{-4}$ & $1\cdot 10^{-5}$ \\
		$\surfT$ & $3\cdot 10^{-3}$ & $3\cdot 10^{-2}$ & $3\cdot 10^{-2}$ & $3\cdot 10^{-2}$  
		\label{tab:CWparamLaStudy}
	\end{tabular}
\end{table}
\end{center}
Note that both phases exhibit the same physical properties. All physical settings are computed on three meshes with equidistant mesh sizes of $\lambda / \grdSz = \{ 8, 16, 32 \}$. The time step sizes on each mesh are set according to the capillary time step restriction \eqref{eq:capillaryTimestep} with a spatial discretization of order $\Pdeg = 2$.

\begin{figure}
	\centering
	\begin{subfigure}{.45\textwidth}	
        \includegraphics[width=.9\linewidth]{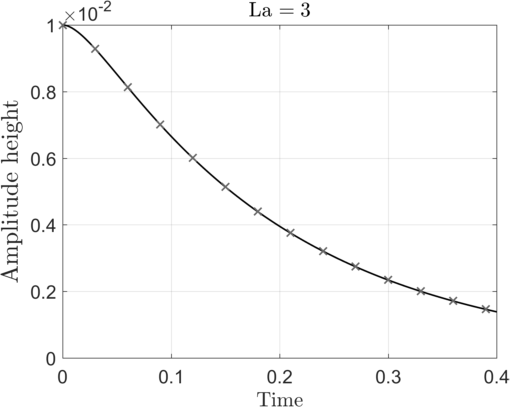}
    \end{subfigure}	
	\begin{subfigure}{.45\textwidth}	
        \includegraphics[width=.9\linewidth]{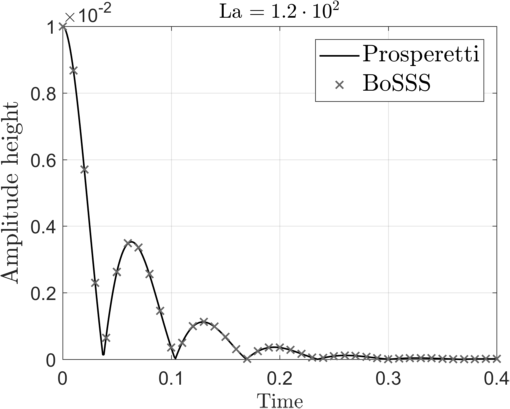}
    \end{subfigure}	
    \\
   	\vspace{10pt}
    \begin{subfigure}{.45\textwidth}	
        \includegraphics[width=.9\linewidth]{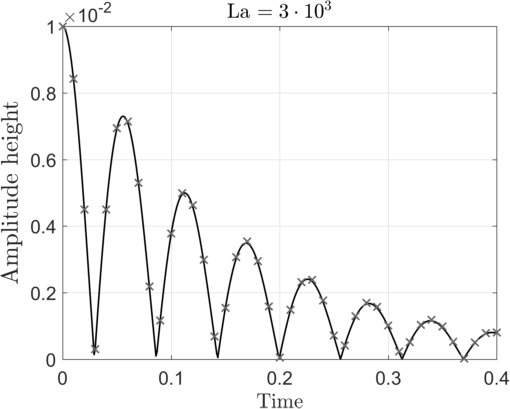}
    \end{subfigure}	
    \begin{subfigure}{.45\textwidth}	
        \includegraphics[width=.9\linewidth]{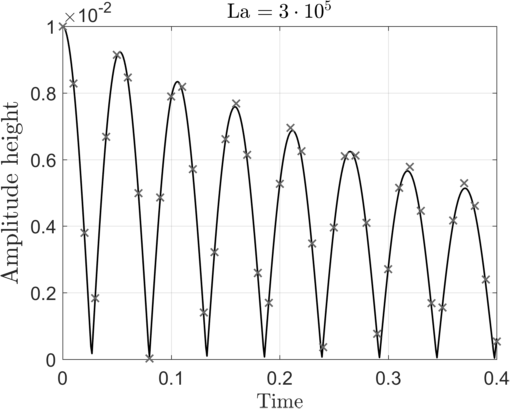}
    \end{subfigure}
	\caption{Comparison of the amplitude height $a$ over time $t$ between the numerical solution of BoSSS (marker set every 100th data point) and the analytical solution\cite{prosperetti_motion_1981} (solid) for a capillary wave with $\laplace = \{ 3, 1.2\cdot 10^2, 3\cdot 10^3, 3\cdot 10^5 \}$.}
	\label{fig:CW_setupStudy}
\end{figure}

In Figure \ref{fig:CW_setupStudy} the numerical solutions on the finest grid are depicted and compared to the analytical reference solution\cite{prosperetti_motion_1981}. The reference solution is computed using MATLAB. Throughout the study all numerical results show very good agreement to the analytical solution. However, one should note that the agreement of the damping rate for the most dynamic case $\textrm{La} = 3 \cdot 10^{5}$ is slightly underestimated resulting in higher amplitude heights for the last peaks. The frequency is still in very good agreement.  

In Table \ref{tab:CW_errors} the $l_2$-errors and the corresponding EOC value for each study are presented. Looking at the finest solutions of each study one should note the increasing error norm for increasing Laplace number. This is additionally attended by decreasing EOC values.  
\begin{center}
	\begin{table}
		\caption{$l_2$-error norms for the capillary wave studies}
		\centering
		\begin{tabular}{c|cccc}
			$\laplace$ & $\lambda / \grdSz = 8$ & $\lambda / \grdSz = 16$ & $\lambda / \grdSz = 32$ & EOC \\
			\hline
			$3$ & $1.067 \cdot 10^-3$ & $2.567 \cdot 10^-2$ &  $3.962 \cdot 10^-3$ & $2.268$ \\
			$1.2 \cdot 10^2$ & $2.954 \cdot 10^-1$ & $1.051 \cdot 10^-2$ & $3.487 \cdot 10^-3$ &  $3.167$ \\
			$3 \cdot 10^3$ & $8.314 \cdot 10^-2$ & $2.700 \cdot 10^-2$ & $1.021 \cdot 10^-2$ &  $1.495$ \\
			$3 \cdot 10^5$ & $1.713 \cdot 10^-1$ & $7.827 \cdot 10^-2$ & $2.084 \cdot 10^-2$ & $1.501$ \\
			$3 \cdot 10^5$ & $1.740 \cdot 10^-1$ & $4.044 \cdot 10^-2$ & $1.288 \cdot 10^-2$ & $1.856$ \\
			(${\dt}_{\textrm{fix}} = 4\cdot10^{-5}$) & & & & \\
			$3 \cdot 10^5$ & $5.068 \cdot 10^-2$ & $8.168 \cdot 10^-3$ & $7.789 \cdot 10^-3$ & $1.336$ \\
			($\Pdeg = 3$, ${\dt}_{\textrm{fix}} = 2\cdot10^{-5}$) & & & & 
			\label{tab:CW_errors}
		\end{tabular}
	\end{table}
\end{center}

Two more studies for the most dynamic case ($\textrm{La} = 3 \cdot 10^5$) are given in Table \ref{tab:CW_errors}. The first one is a mesh study performed with a fixed time step size for all simulations that corresponds to a capillary time step size with an grid size of $\lambda / \grdSz = 64$. In this case the error norms and the EOC value are slightly improved. For the second study the polynomial degree is set to $\Pdeg = 3$ and again the time step size is fixed. As expected the error norms for $\Pdeg = 3$ are smaller compared to $\Pdeg = 2$, but the study shows a poorer EOC value. This may result from the theoretically insufficient integration scheme, i.e third order BDF, compared to the spatial approximation in context of the moving interface approach\cite{kummer_time_2018}. Furthermore, one should note that quite small time steps are needed in order to obtain the desired convergence rates. 

\subsection{Droplet simulations}
\label{sec:droplet}

In this section two settings for the transient simulation of a droplet are presented. The first test case considers a droplet in its equilibrium state, i.e. circular shape. For the second one the droplet is initially deformed into an ellipsoidal shape.   

\subsubsection{Droplet in equilibrium}

The equilibrium state, i.e. $\veloc = 0$, of a droplet in $\mathbb{R}^2$ is described by an circular shape, where the pressure inside a droplet (phase $\domainA$) of radius $r$ is given by the Young-Laplace equation with 
\begin{equation}
    \indA{\press} = \indB{\press} + \frac{\surfT}{r}.
\end{equation}
However, due to numerical inaccuracies regarding the surface tension force the computed solutions do not provide a zero-velocity field. In order to quantify the discretization error Smolianksi \cite{smolianski_numerical_2001} proposed the following setup. The droplet with radius $r = 0.25$ is set in the middle of the computational domain of $\domain = [-0.5, 0.5] \times [-0.5, 0.5]$, where all boundaries describe a wall boundary condition. Both phases share the same physical properties with: density $\dens = 10^4$, dynamic viscosity $\visc = 1$ and surface tension coefficient $\surfT = 1$. This results in a Laplace number of $\laplace = \surfT \dens 2 r / \visc^2 = 5 \cdot 10^3$. We investigate a mesh study with the following mesh sizes $1/h = \{ 20, 40, 60, 80\}$. All simulations are performed until $t = 125$ with a fixed time step size of $\dt = 0.01$. Setting the polynomial degree to $\Pdeg = 2$ the capillary time step restriction \eqref{eq:capillaryTimestep} is satisfied for all runs.

\begin{figure}
	\centering
	\begin{subfigure}{.45\textwidth}	
        \includegraphics[width=.9\linewidth]{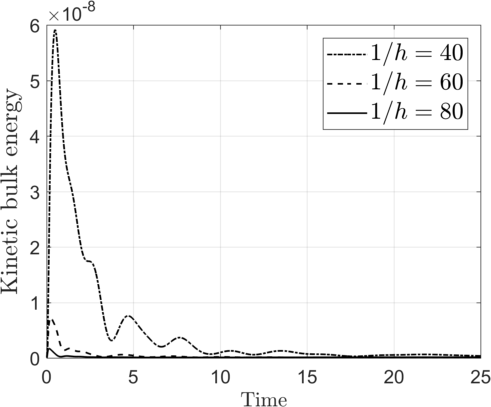}
    \end{subfigure}	
	\begin{subfigure}{.45\textwidth}	
        \includegraphics[width=.9\linewidth]{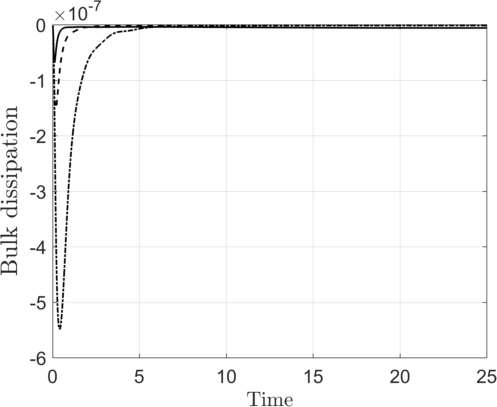}
    \end{subfigure}	
    \begin{subfigure}{.45\textwidth}	
        \includegraphics[width=.9\linewidth]{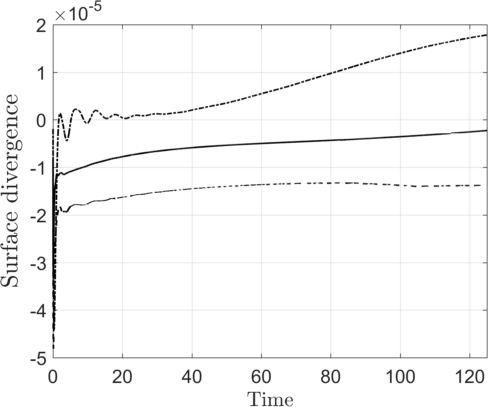}
    \end{subfigure}
	\caption{The temporal evolution of the $L^2$-error norms of the kinetic bulk energy $\frac{1}{2}\LtwoNorm{\dens\veloc \cdot \veloc}$ (left), the bulk dissipation $\LtwoNorm{\textrm{tr}(\tensr{D}(\veloc)^2)}$ (middle) and the surface divergence $\norm{\divI{\velocI}}_{L^2,\interface}$ (right) for a droplet in equilibrium state. The transient simulations are performed on meshes with $1 /\grdSz = \{ 40, 60, 80\}$.}
	\label{fig:SD_energies}
\end{figure}

In Table \ref{tab:SD_LaplaceError} the $L^2$-error norm of the spurious velocities $\norm{\veloc}_{L^2}$ and the error against the exact Young-Laplace solution $\norm{\press - \press_\textrm{exact}}_{L^2}$ are given. The norm is computed in the bulk $\bulk$ on the terminal time step at $t = 125$.
\begin{center}
\begin{table}
\caption{$L^2$-error norms of spurious velocities and Laplace-Young equation for the terminal time step at $t = 125$.}
\centering
	\begin{tabular}{c|cccc}
		$1/\grdSz$ & 20 & 40 & 60 & 80 \\
		\hline
		$\norm{\veloc}_{L^2}$ & $1.68\cdot 10^{-5}$ & $2.95\cdot 10^{-7}$ & $2.60\cdot 10^{-7}$ & $1.34\cdot 10^{-7}$ \\
		$\norm{\press - \press_\textrm{exact}}_{L^2}$ & $1.03\cdot 10^{-2}$ & $5.20\cdot 10^{-4}$ & $5.68\cdot 10^{-4}$ & $4.15\cdot 10^{-4}$
		\label{tab:SD_LaplaceError}
	\end{tabular}
\end{table}
\end{center}
One should note that the error norms seem to converge already on the second mesh. This may result from the long simulation time and fixed time step size. However, we take a closer look on some energy related properties. In Figure \ref{fig:SD_energies} the kinetic energy in the bulk $\frac{1}{2}\LtwoNorm{\dens \veloc \cdot \veloc}$ and the corresponding discrete dissipation term $\LtwoNorm{\textrm{tr}(\tensr{D}(\veloc)^2)}$ are depicted for $t = [0, 25]$. Furthermore, the surface divergence $\norm{\divI{\velocI}}_{L^2,\interface}$ is given for $t = [0, 125]$.  
One observes that the numerical solutions on the coarser meshes exhibit a larger initial deviation from the zero-velocity field and that the decay of such spurious velocities takes a longer period of time. The same behavior is shown for the bulk dissipation. One should note that the term is always negative demonstrating the energetic stability of the presented discretization. However, taking the surface divergence, i.e. change rate of the interface area, into account the coarser meshes exhibit an enlarging or diminishing interface area until the end. The finest solution seem to produce a stable interface.  

\subsubsection{Oscillating droplet}

In this section following the work of Hysing \cite{hysing_new_2006} we consider a droplet initialized in an ellipsoidal shape with 
\begin{equation}
	\levset = \frac{x^2}{a^2} + \frac{y^2}{b^2} - 1,
	\label{eq:ellipsoid}
\end{equation}
where the semi-axes are set to $a = 1.25r$ and $b = 0.8r$. The corresponding equilibrium shape is described by an circle with radius $r = 0.25$. The computational domain and the physical properties are same as for the equilibrium test case, except for the surface tension coefficient with $\surfT = 0.1$. This results in a smaller Laplace number of $\textrm{La} \approx 500$. A mesh study on the following mesh sizes $1/h = \{ 10, 20, 40, 60, 80\}$ is performed and the simulations run until $t = 100$ with a fixed time step size of $\dt = 0.5$. Setting the polynomial degree to $\Pdeg = 2$, the capillary time step restriction is already exceeded on the second mesh.

\begin{figure}
	\centering
	\begin{subfigure}{.45\textwidth}	
        \includegraphics[width=.9\linewidth]{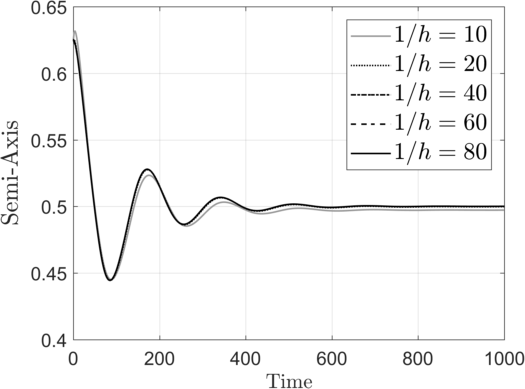}
    \end{subfigure}	
	\begin{subfigure}{.45\textwidth}	
        \includegraphics[width=.9\linewidth]{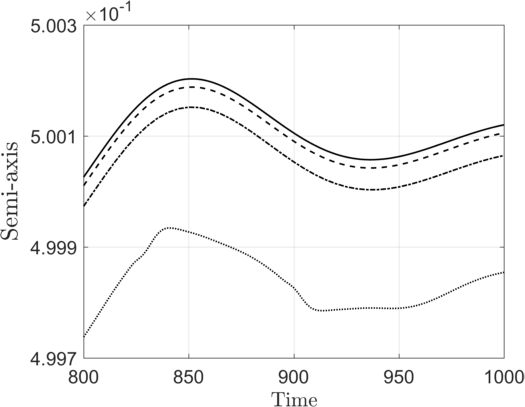}
    \end{subfigure}	
    \begin{subfigure}{.45\textwidth}	
        \includegraphics[width=.9\linewidth]{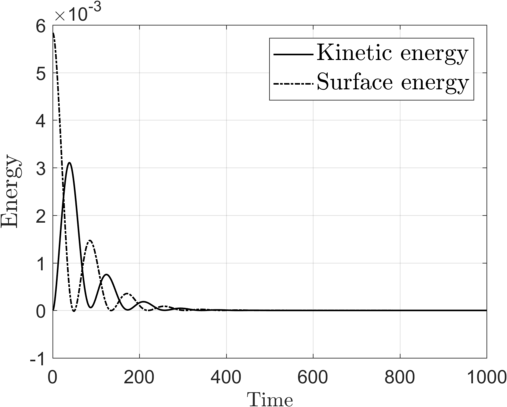}
    \end{subfigure}
	\caption{The temporal evolution of the semi-axis in $x$-direction of an oscillating droplet for $t = [0, 1000]$ (left) and $t = [800, 1000]$ (middle). On the right the temporal evolution of the kinetic bulk energy and surface energy for the finest solution with $1/\grdSz = 80$.}
	\label{fig:OD_meshStudy}
\end{figure}

In Figure \ref{fig:OD_meshStudy} on the left the temporal evolution of the semi-axis in $x$-direction $a_x$ is plotted for the whole mesh study. It is remarkable that all simulations show a stable oscillating behavior until the end, see Figure \ref{fig:OD_meshStudy} in the middle. Note that on the finest mesh the time step size is 10 times larger than the capillary time step restriction \eqref{eq:capillaryTimestep}.

Taking a closer look on the finest solution, the temporal evolution of the kinetic bulk energy and the surface energy $E_{\surfT} = \surfT \abs{\interface}$ with $\abs{\interface}$ describing the interface area is plotted in Figure \ref{fig:OD_meshStudy} on the right. Note that the plotted surface energy is subtracted by the minimal surface, i.e. $E_{\surfT, \textrm{min}} = \surfT 2 \pi r$. Thus, the exchange of kinetic and surface energy is observable in one plot. The maximum values of the surface energy, i.e. largest extend along one semi-axis, corresponds directly to the minimum value of the kinetic energy. Between these states the kinetic energy reaches its maximum always slightly before the minimum value of the surface energy. Due to the small Laplace number the oscillations are strongly damped by viscous dissipation in the bulk.

\subsection{Rising bubble benchmark}
\label{sec:RB}

\begin{figure}
	\centering
	\def\svgwidth{120pt}
	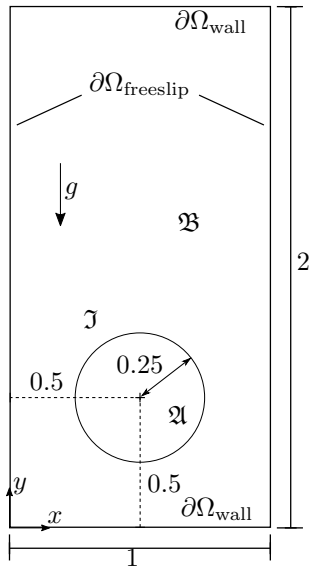
	\caption{Numerical setup for the rising bubble benchmark test case.}
	\label{fig:RBsetup}
\end{figure}

In this last section we compare our results against the rising bubble benchmark test case established by Hysing et al.\cite{hysing_quantitative_2009}. The initial setup is depicted in Figure \ref{fig:RBsetup}, where a circular bubble (radius $r = 0.25$) is positioned at $(0.5, 0.5)$ in the domain $\domain = [0,1] \times [0,2]$. 
The lower and upper boundary are imposed by the wall boundary condition. The boundaries at the sides describe a free-slip boundary condition given by 
\begin{equation}
	\veloc \cdot \normalB = 0, \quad \tangent_{\boundary} \cdot \left( \grad{\veloc} + \gradT{\veloc} \right) \normalB = 0,
	\label{eq:freeslip}
\end{equation}
where $\tangent_{\boundary}$ denotes the tangent vector on the boundary $\partial \domain$. The discretization in context of the presented variational problem \eqref{eq:variationalForm_NSE} is given in the Appendix \ref{app:FreeSlipBC}. The driving force for the rise of the bubble ( $\indA{\dens} < \indB{\dens}$) is the gravity force $\vectr{g} = -g \vectr{e}_y$ oriented in negative $y$-direction. 

For this initial setting two test cases are defined with different physical properties resulting in an ellipsoidal terminal shape for the first test case and a dimpled cap with filaments for the second test case. The corresponding physical parameters are given in Table \ref{tab:physParamRB}. For both test cases the simulation is performed for 3 time units.
\begin{center}
\begin{table}
\caption{Physical parameters for both test cases of the rising bubble benchmark.}
\centering
	\begin{tabular}{c|cccccc}
		Test case & $\indA{\dens}$ & $\indB{\dens}$ & $\indA{\visc}$ & $\indB{\visc}$ & $\surfT$ & $g$\\
		\hline 
		1 & 100 & 1000 & 1 & 10 & 24.5 & 0.98 \\
		2 & 1 & 1000 & 0.1 & 10 & 1.96 & 0.98
		\label{tab:physParamRB}
	\end{tabular}
\end{table}
\end{center}

\paragraph{Benchmark groups and comparing quantities} 

In Hysing et al.\cite{hysing_quantitative_2009} an extensive data-set is provided by the following three research groups: TP2D (Transport Phenomena in 2D), FreeLIFE (Free-Surface Library of Finite Element) and MoonMD (Mathematics and object-oriented Numerics in Magdeburg). For details on the methodology of each solver the reader is referred to the benchmark paper\cite{hysing_quantitative_2009}. The plotted data of the benchmark groups are taken from \url{http://www.featflow.de/en/benchmarks/cfdbenchmarking/bubble/bubble_reference.html}. Furthermore, we compare in the following some results to the work of Heimann et al.\cite{heimann_unfitted_2013}, where the unfitted DG method with a piece-wise linear approximation of the interface is employed. 

In order to compare the different numerical methods three scalar measures for the temporal evolution of the rising bubble are introduced\cite{hysing_quantitative_2009}. First, the center of mass is considered given as 
\begin{equation}
	\vectr{x}_{\textrm{c}} = \frac{\int_{\domainA} \vectr{x} \dif{V}}{\int_{\domainA} 1 \dif{V}}.
	\label{eq:centerOfMass}
\end{equation}
The second measure is denoted as circularity and defined by 
\begin{equation}
	\cancel{c} = \frac{\textrm{perimeter of area-equivalent circle}}{\textrm{perimeter of bubble}}.
	\label{eq:circularity}
\end{equation}
The circularity is $\cancel{c} = 1$ for a perfectly circular bubble and gets smaller for more deformed bubble shapes. The third measure is the mean bubble velocity 
\begin{equation}
	\veloc_\textrm{c} = \frac{\int_{\domainA} \veloc \dif{V}}{\int_{\domainA} 1 \dif{V}},
	\label{eq:meanBubbleVeloc}
\end{equation}
where the component in $y$-direction is denoted as the rise velocity $V_c$. The error quantification is done via the $l_2$-error norm \eqref{eq:ScalarErrorNorms2} and the corresponding ROC values \eqref{eq:ConvergenceRate}. Note that the in the original benchmark paper furthermore $l_1$ and $l_\infty$-error norms are considered. For completeness the definitions and corresponding results are found in the Appendix \ref{app:risingBubble}. Besides the three scalar quantities the terminal shape of the bubble at $t = 3$ is compared.

\subsubsection{Test case 1 - Ellipsoidal shape}

\begin{figure}
	\centering
	\pgfplotsset{width=0.4\textwidth}
	\subfloat{	
		\begin{tikzpicture}
		\begin{axis}[
		enlargelimits=false,
		grid=none,
		trim axis left,
		trim axis right,
		axis equal image,
		axis on top,
		xtick distance=0.5,
		xlabel=$x$, 
		ylabel=$y$,
		]
		\addplot graphics [xmin=0, xmax=1, ymin=0, ymax=2]
		{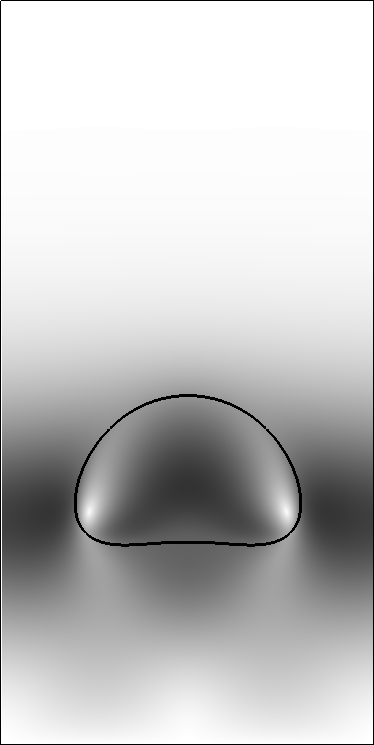};
		\end{axis}
		\end{tikzpicture}
	}
	\subfloat{	
		\begin{tikzpicture}
		\begin{axis}[
		enlargelimits=false,
		grid=none,
		trim axis left,
		trim axis right,
		axis equal image,
		axis on top,
		xtick distance=0.5,
		xlabel=$x$, 
		ytick=\empty,
		%
		]
		\addplot graphics [xmin=0, xmax=1, ymin=0, ymax=2]
		{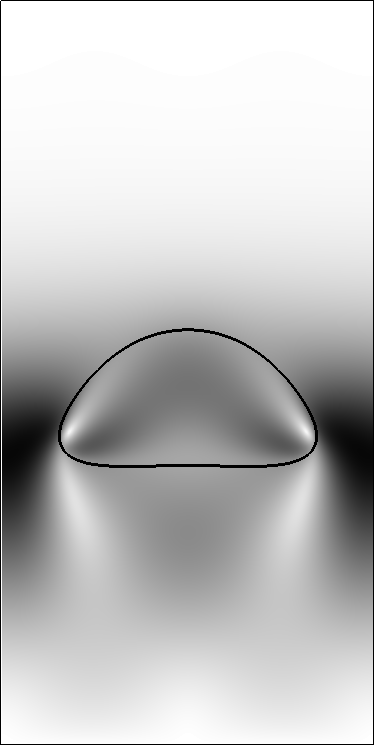};
		\end{axis}
		\end{tikzpicture}
	}
	\subfloat{	
		\begin{tikzpicture}
		\begin{axis}[
		enlargelimits=false,
		grid=none,
		trim axis left,
		trim axis right,
		axis equal image,
		axis on top,
		xtick distance=0.5,
		xlabel=$x$, 
		ytick=\empty,
		%
		colorbar,
		colormap={}{ gray(0cm)=(1); gray(1cm)=(0);},
		colorbar style={
			align=right,
			ylabel = $\abs{\veloc}$,
			/pgf/number format/fixed,
			ytick={0,0.25,0.5},
			tick align=outside,
			tick pos=right
		},
		point meta min=0,
		point meta max=0.5,
		]
		\addplot graphics [xmin=0, xmax=1, ymin=0, ymax=2]
		{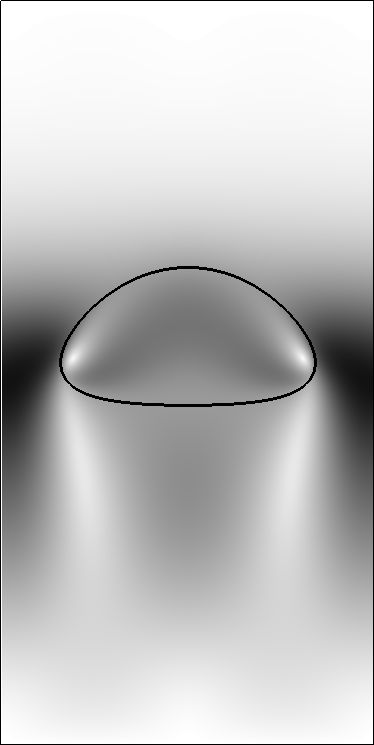};
		\end{axis}
		\end{tikzpicture}
	}
	\caption{Temporal evolution of the rising bubble benchmark test case 1 for times $t = \{ 1.2, 2.1, 3.0\}$. The shown results are computed with a polynomial degree of $\Pdeg = 3$ on a mesh with $1/\grdSz = 60$. The plotted field describes the magnitude of the velocity vector field $\abs{\veloc}$.}
	\label{fig:RBevolution_tc1}
\end{figure}

For the first test case two mesh studies were performed for polynomial degrees of $\Pdeg = \{ 2, 3 \}$ on meshes with equidistant cell sizes of $1/\grdSz = \{ 10, 20, 40, 60, 80 \}$. The corresponding number of background cells, total number of DOF and the number of time steps are given in Table \ref{tab:RBtc1NDOF}. The total number of DOF indicates the maximum number encountered during the simulation. Since the number of cut-cells is depending on the interface form, the number of double DOF is changing according to the number of cut-cells. The time step sizes are chosen for each setting according to the capillary time step restriction \eqref{eq:capillaryTimestep}.\\
\begin{center}
\begin{table}
\caption{Number of background cells (NEL), total number of DOF (NDOF) and number of time steps (NTS) for the rising bubble benchmark test case 1 for polynomial degrees $\Pdeg = \{2, 3\}$.}
\centering
	\begin{tabular}{c|ccccc}
	    & & \multicolumn{2}{c}{$\Pdeg = 2$ } &  \multicolumn{2}{c}{$\Pdeg = 3$ }\\
		$1/\grdSz$ & NEL & NDOF & NTS & NDOF & NTS \\
		\hline 
		10 & 200 & - & - & 5564 & 500 \\
		20 & 800 & 12390 & 600 & 21476 & 1000 \\
		40 & 3200 & 48810 & 2000 & 84604 & 3000 \\
		60 & 7200 & 109230 & 3000 & 189332 & 5000 \\
		80 & 12800 & 192650 & 5000 & -  & - 
		\label{tab:RBtc1NDOF}
	\end{tabular}
\end{table}
\end{center}
The temporal evolution of the bubble shape at times $t = \{ 1.2, 2.1, 3.0\}$ is shown in Figure \ref{fig:RBevolution_tc1}. Here, the magnitude of the corresponding velocity field is plotted for the finest solution of $\Pdeg = 3$. Comparing the terminal shape at $t=3$ with the respective finest solution of benchmark groups (Figure \ref{app:RBtc1_terminalShape} in Appendix \ref{app:risingBubble}), there are no significant differences visible throughout all methods.\\

The $l_2$-error norms for the benchmark quantities against the finest solution and the corresponding ROC values on each refinement level are given in Table \ref{tab:RB_convStudy_errorNorms2}. The ROC values of all benchmark quantities for $\Pdeg = 2$ show a high-order convergence up to 2.5, whereas the study for $\Pdeg = 3$ shows poorer results with ROC below 2. This is again a result regarding the moving interface approach (see Section \ref{sec:CW}), since the temporal discretization with a third order BDF does not match the theoretical requirement for a time integration scheme of at least $2\Pdeg$. Thus, for $\Pdeg = 3$ we cannot expect the higher-order convergence rates. However, the higher polynomial order still provides smaller error norms on the same mesh sizes. Note that a reinitialization is performed every 50th time step on the coarsest mesh for $\Pdeg = 3$ in order to stabilize the interface evolution. This additional operation may lead to the lower error norm compared to the next finer mesh resulting in a negative ROC value. Note that the error norms $l_1$ and $l_\infty$ with the corresponding ROC values are found in the Appendix \ref{app:risingBubble} Tables \ref{app:RB_convStudy_errorNorms1} and \ref{app:RB_convStudy_errorNormsInfty}. The above statements are also valid for both other error norms.\\ 
\begin{center}
\begin{table}
\caption{$l_2$-error norms against the finest solution and $\textrm{ROC}$ for the mesh studies ($\Pdeg = \{ 2, 3\}$) of the rising bubble benchmark test case 1.}
\centering
	\begin{tabular}{cc|cc|cc|cc}
		&& \multicolumn{2}{c}{Center of mass} & \multicolumn{2}{c}{Rise velocity} & \multicolumn{2}{c}{Circularity} \\
		&$1/\grdSz$ & $\norm{e}_2$ & $\textrm{ROC}$ & $\norm{e}_2$ & $\textrm{ROC}$ & $\norm{e}_2$ & $\textrm{ROC}$ \\
		\hline 
		\hline
		\multirow{3}{*}{\rotatebox[origin=c]{90}{$\Pdeg = 2$}} & 20 & $4.92 \cdot 10^{-3}$ & - & $1.02 \cdot 10^{-2}$ & - & $8.77 \cdot 10^{-4}$ & - \\
		& 40 & $1.27 \cdot 10^{-3}$ & 1.95 & $2.95 \cdot 10^{-3}$ & 1.79 & $5.49 \cdot 10^{-4}$ & 0.68 \\
		& 60 & $4.53 \cdot 10^{-4}$ & 2.54 & $1.05 \cdot 10^{-3}$ & 2.55 & $1.92 \cdot 10^{-4}$ & 2.59 \\
		\hline
		\multirow{3}{*}{\rotatebox[origin=c]{90}{$\Pdeg = 3$}} & 10 & $1.44 \cdot 10^{-2}$ & - & $1.56 \cdot 10^{-2}$ & - & $9.58 \cdot 10^{-4}$ & -  \\
		& 20 & $2.23 \cdot 10^{-3}$ & 2.69 & $5.96 \cdot 10^{-3}$ & 1.39 & $1.65 \cdot 10^{-3}$ & -0.78 \\
		& 40 & $6.85 \cdot 10^{-4}$ & 1.70 & $1.69 \cdot 10^{-3}$ & 1.82 & $4.67 \cdot 10^{-4}$ & 1.04
		\label{tab:RB_convStudy_errorNorms2}
	\end{tabular}
\end{table}
\end{center}

\begin{figure}
	\centering
	\begin{subfigure}{.45\textwidth}	
        \includegraphics[width=.9\linewidth]{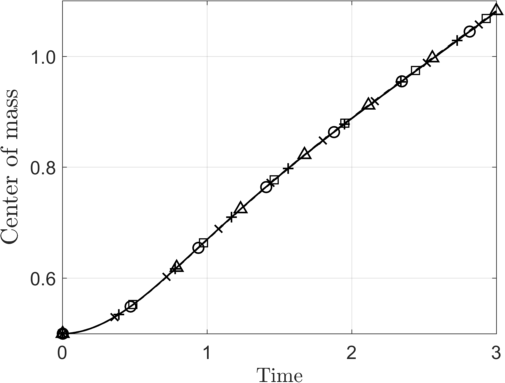}
    \end{subfigure}	
	\begin{subfigure}{.45\textwidth}	
        \includegraphics[width=.9\linewidth]{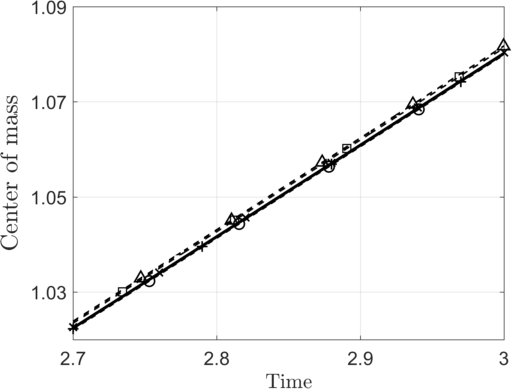}
    \end{subfigure}	
    \par\bigskip
	\begin{subfigure}{.45\textwidth}	
        \includegraphics[width=.9\linewidth]{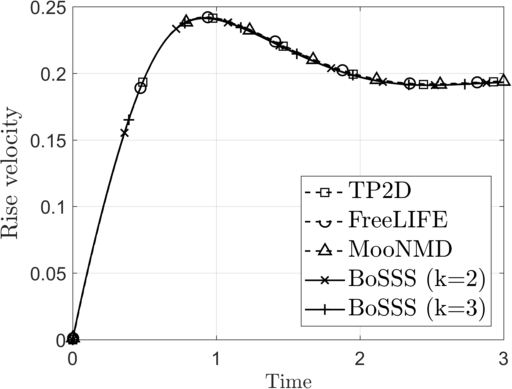}
    \end{subfigure}	
	\begin{subfigure}{.45\textwidth}	
        \includegraphics[width=.9\linewidth]{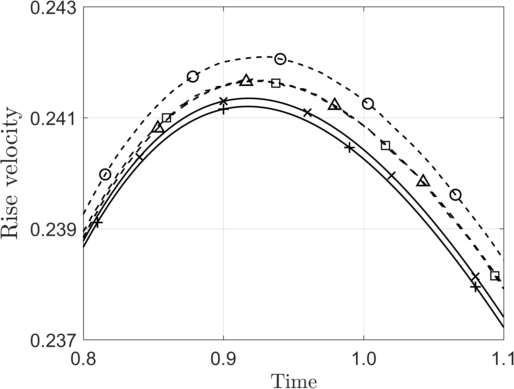}
    \end{subfigure}	
    \par\bigskip
	\begin{subfigure}{.45\textwidth}	
        \includegraphics[width=.9\linewidth]{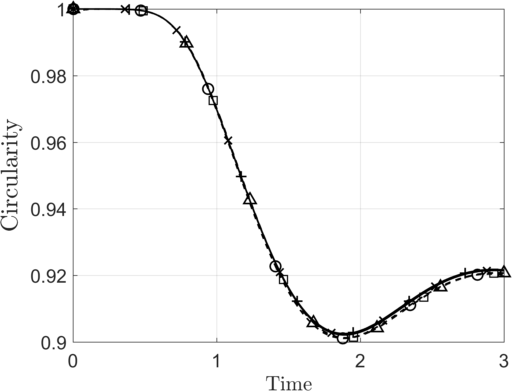}
    \end{subfigure}	
	\begin{subfigure}{.45\textwidth}	
        \includegraphics[width=.9\linewidth]{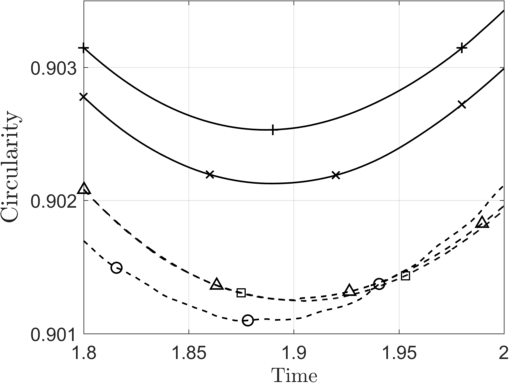}
    \end{subfigure}	
	\caption{Comparison to the benchmark groups for the temporal evolution of the $y$-component of the center of mass $y_c$ (top row), of the rise velocity $V_c$ (middle row) and the circularity $\cancel{c}$ (bottom row). In the left column the overall solutions are depicted and the right column provides a zoom for specific time intervals. The plotted solutions show the finest solution of each group for test case 1.}
	\label{fig:RBtc1_benchmarkQuantities}
\end{figure}

The results of the benchmark quantities on the finest mesh for each study are compared in the following against the finest solutions of the benchmark groups, see Figure \ref{fig:RBtc1_benchmarkQuantities}. In the top row of Figure \ref{fig:RBtc1_benchmarkQuantities} the temporal evolution of the center of mass $y_c(t)$ is plotted. Overall, all numerical solutions agree excellently on the evolution, only a zoom for $t = [2.75, 3.0]$ exhibits two distinct trends. Both our numerical solutions coincide very well with FreeLIFE. Likewise, the other two benchmark groups coincide, but are constantly above ours. The terminal rise height at $t = 3$ for both studies ($\Pdeg = \{2,3\}$) is given in Appendix \ref{app:FreeSlipBC}, see Table \ref{app:RBtc1_distinctValues}. Besides our results, the reader finds the results of the benchmark groups and of Heimann et al.\cite{heimann_unfitted_2013} for a direct comparison. All results are in the range of $y_c(t=3) = [1.0799, 1.0817]$.\\

The temporal evolution of the rise velocity $V_c(t)$ depicted in Figure \ref{fig:RBtc1_benchmarkQuantities} in the middle row shows a very agreement between all groups. A closer look at the maximum velocity at around $t = 0.9$ displays that our numerical solutions slightly underestimate the peak velocity compared to the benchmark groups. However, one should note that FreeLIFE overestimates both other benchmark groups in the same range, i.e. the deviation is around $ 4\cdot10^{-4}$. Furthermore, one observes that higher peak values are obtained later in time, see Table \ref{app:RBtc1_distinctValues} in the Appendix. Looking at Heimann et al.\cite{heimann_unfitted_2013}, their result also tends to underestimate the benchmark groups.\\ 

Considering the circularity (Figure \ref{fig:RBtc1_benchmarkQuantities} in bottom row) the overall qualitative agreement is still very good, but our solutions clearly overestimate the minimum value at around $t = 1.9$. Zooming in, the results of TP2D and MooNMD indistinguishably coincide. FreeLIFE reaches its minimum value slightly before the other groups. Note that our results are qualitatively closer to the solution of FreeLIFE exhibiting a more gentle slope in the beginning and getting steeper to the end. A reason for the overestimation may be the mass production during the simulation which is, however, comparably small for the finest solution, i.e. the relative mass error is around $10^{-3}$. The minimum circularity and the corresponding time are given in the Appendix in Table \ref{app:RBtc1_distinctValues}.

\subsubsection{Test case 2 - Dimpled cap with filaments}

\begin{figure}
	\centering
	\begin{minipage}{0.4\textwidth}
	\centering
	\pgfplotsset{width=\textwidth}
		\begin{tikzpicture}
		\begin{axis}[
		enlargelimits=false,
		grid=none,
		trim axis left,
		trim axis right,
		axis equal image,
		axis on top,
		xtick distance=0.5,
		xlabel=$x$, 
		ylabel=$y$,
		]
		\addplot graphics [xmin=0, xmax=1, ymin=0, ymax=2]
		{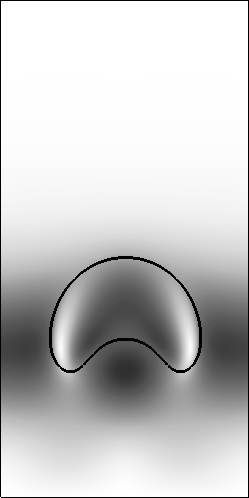};
		\end{axis}
		\end{tikzpicture}
		\begin{tikzpicture}
		\begin{axis}[
		enlargelimits=false,
		grid=none,
		trim axis left,
		trim axis right,
		axis equal image,
		axis on top,
		xtick distance=0.5,
		xlabel=$x$, 
		ytick=\empty,
		%
		]
		\addplot graphics [xmin=0, xmax=1, ymin=0, ymax=2]
		{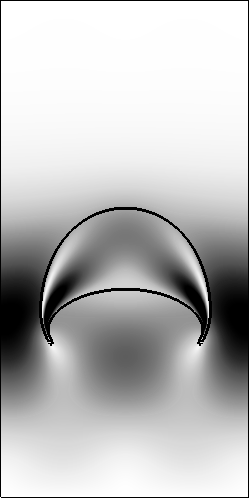};
		\end{axis}
		\end{tikzpicture}
		\begin{tikzpicture}
		\begin{axis}[
		enlargelimits=false,
		grid=none,
		trim axis left,
		trim axis right,
		axis equal image,
		axis on top,
		xtick distance=0.5,
		%
		colorbar,
		colormap={}{ gray(0cm)=(1); gray(1cm)=(0);},
		colorbar style={
			align=right,
			ylabel = $\abs{\veloc}$,
			/pgf/number format/fixed,
			ytick={0, 0.25, 0.5},
			tick align=outside,
			tick pos=right
		},
		point meta min=0,
		point meta max=0.5,
		]
		\addplot graphics [xmin=0, xmax=1, ymin=0, ymax=2]
		{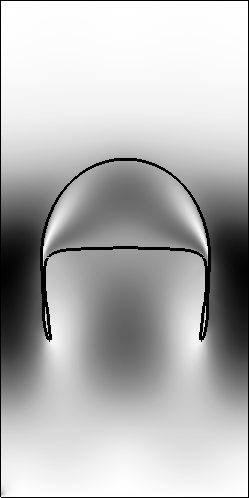};
		\end{axis}
		\end{tikzpicture}
	\end{minipage}
	\begin{minipage}{0.55\textwidth}
	\centering
	\pgfplotsset{width=1.2\textwidth}
		\begin{tikzpicture}
		\begin{axis}[
		enlargelimits=false,
		grid=none,
		trim axis left,
		trim axis right,
		axis equal image,
		axis on top,
		xtick distance=0.1,
		xlabel=$x$, 
		ylabel=$y$,
		]
		\addplot graphics [xmin=0.7, xmax=0.9, ymin=0.55, ymax=1.05]
		{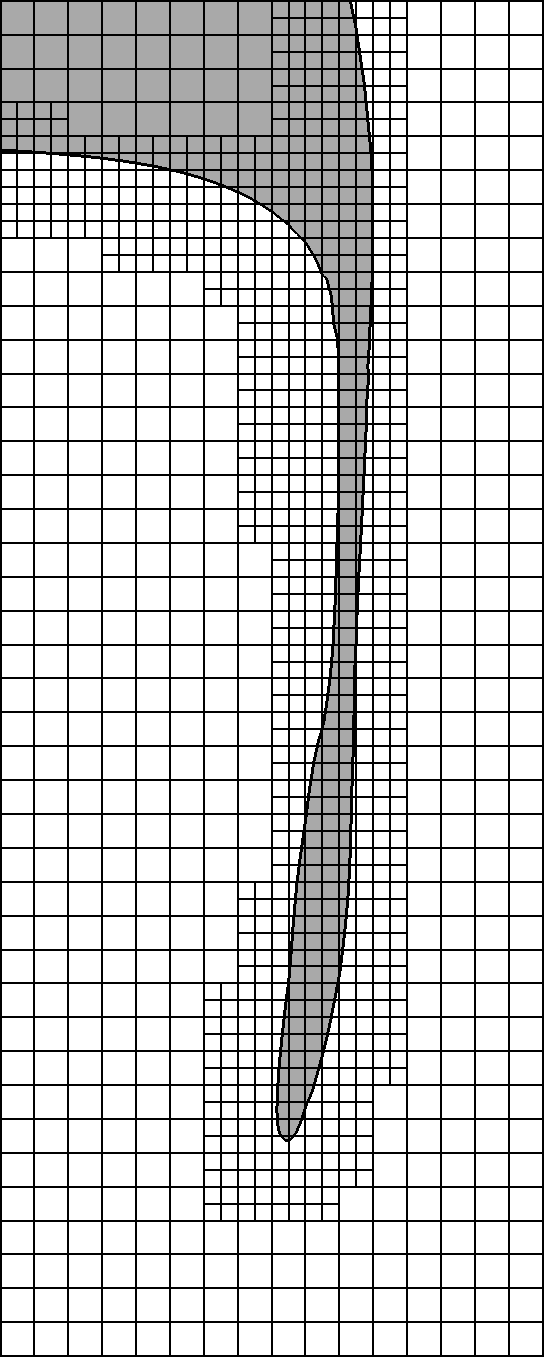};
		\end{axis}
		\end{tikzpicture}
		\begin{tikzpicture}
		\begin{axis}[
		enlargelimits=false,
		grid=none,
		trim axis left,
		trim axis right,
		axis equal image,
		axis on top,
		xtick distance=0.1,
		xlabel=$x$, 
		ytick=\empty,
		]
		\addplot graphics [xmin=0.7, xmax=0.9, ymin=0.55, ymax=1.05]
		{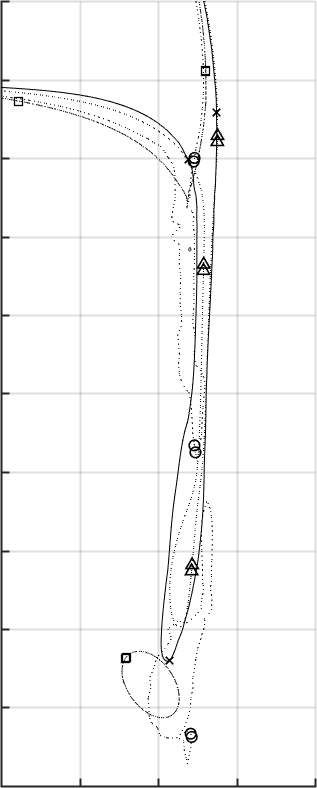};
		\end{axis}
		\end{tikzpicture}
	\end{minipage}
	\caption{Left: Temporal evolution of the rising bubble benchmark test case 2 for times $t = \{ 1.2, 2.1, 3\}$. The shown results are computed with a polynomial degree of $\Pdeg = 2$ on a mesh with $1/\grdSz = 80$ and AMR level 1 on the narrow band. The plotted field describes the magnitude of the velocity vector field $\abs{\veloc}$. Middle: Close-up on the right filament showing the actual numerical mesh ($1/\grdSz = 80$) with AMR level 1. Right: Comparison between the benchmark groups. BoSSS (crosses) with $1/\grdSz_\textrm{min} = 160$, TP2D (squares) with $1/\grdSz = 640$, FreeLIFE (circles) with $1/\grdSz = 160$, MooNMD (triangles) with $\textrm{NDOF}_\textrm{int} = 900$.}
	\label{fig:RBevolution_tc2}
\end{figure}

For the second test case a mesh study with $\Pdeg = 2$ is done on meshes with $1/\grdSz = \{ 20, 40, 80\}$. All simulations are performed with an additional constant mesh refinement on the narrow band with level 1, since the bubble shape exhibits stronger deformations for this case. In Figure \ref{fig:RBevolution_tc2} the temporal evolution at times $t = \{ 1.2, 2.1, 3.0\}$ is depicted on the left for the finest solution, where the magnitude of velocity field $\abs{\veloc}$ is plotted. Compared to the first test case (Figure \ref{fig:RBevolution_tc1}) the bubble exhibits a considerably more concave deformation. This evolves into a cap-like shape,
where thin filaments starts to emerge from the bubble and get longer
and thinner close to the bubble. The terminal shape at $t=3$ of the bubble assumes a dimpled cap. A close-up of the right filament is displayed in the middle of Figure \ref{fig:RBevolution_tc2}, where the actual mesh of our simulation is shown. Note that there is a buffer layer of refined cells outside the narrow band, which mark previous cells in the narrow band. The thickness of the filament in gray is determined by the spatial resolution of the adapted mesh. 

A comparison between the benchmark groups is given on the right in Figure \ref{fig:RBevolution_tc2}. Since such filaments are strongly mesh dependent, one cannot expect good agreement in such regions. The result of TP2D ($1/\grdSz = 320$) even exhibit a break up of the filament with an additional satellite droplets. Taking a closer look on the transition between the main bubble and the filament, two different shapes may be characterized. The filament is slightly more extended to the outside of the bubble in our solution and the one from MooNMD. On the other side TP2D and FreeLIFE show a more inward curved shape at the transition region.\\ 

Comparing the benchmark quantities, the agreement for the center of mass between all methods is again very good, see Figure \ref{fig:RBtc2_CenterOfMass} in the Appendix \ref{app:risingBubble}. The temporal evolution is not much changed compared to the first test case, but the deviations to the end are larger, where our numerical result lies in between the benchmark groups. The terminal rise height and the other distinct values for the benchmark groups and Heimann et al.\cite{heimann_unfitted_2013} are again found in the Appendix \ref{app:risingBubble} in Table \ref{app:RBtc2_distinctValues}. 

\begin{figure}
	\centering
	\begin{subfigure}{.45\textwidth}	
        \includegraphics[width=.9\linewidth]{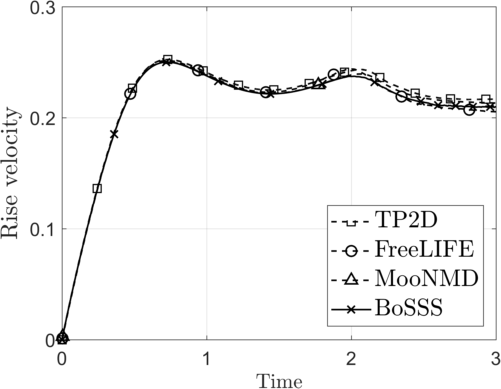}
    \end{subfigure}	
	\begin{subfigure}{.45\textwidth}	
        \includegraphics[width=.9\linewidth]{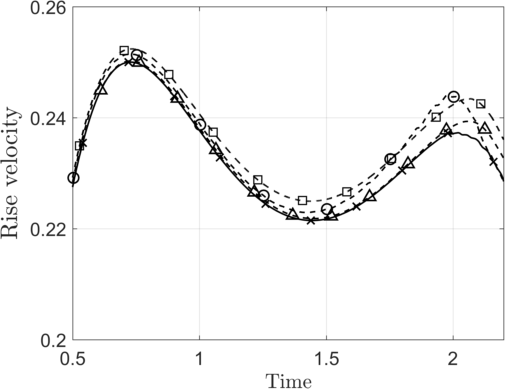}
    \end{subfigure}
    \par\bigskip
    \begin{subfigure}{.45\textwidth}	
        \includegraphics[width=.9\linewidth]{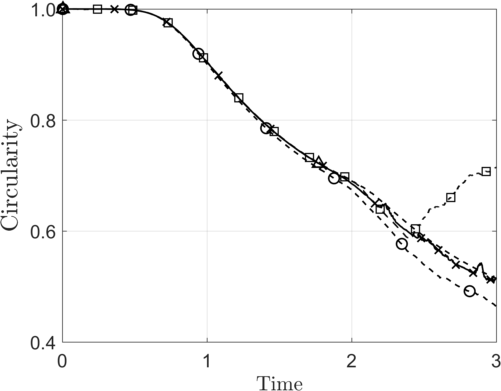}
    \end{subfigure}	
	\begin{subfigure}{.45\textwidth}	
        \includegraphics[width=.9\linewidth]{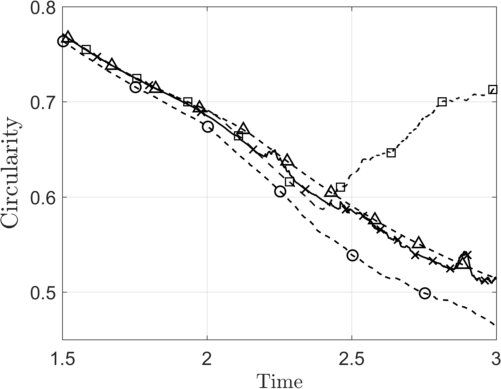}
    \end{subfigure}
	\caption{Comparison to the benchmark groups for the temporal evolution of the rise velocity $V_c$ (top row) and the circularity $\cancel{c}$ (bottom row). In the left column the overall solutions are depicted and the right column provides a zoom for specific time intervals. The plotted solutions show the finest solution of each group for test case 2.}
	\label{fig:RBtc2_benchmarkQuantities}
\end{figure}

The temporal evolution of the rise velocity (Figure \ref{fig:RBtc2_benchmarkQuantities} in top row) exhibits an additional peak around $t = 2$, which is captured by all methods. The second peak is due to the emerging filaments form the main bubble, see Figure \ref{fig:RBevolution_tc2}. Taking a closer look on both peaks, the deviations at the second one are clearly larger again revealing two distinct trends. Our solution agrees very well with MooNMD underpredicting the whole evolution of rise velocity compared to both other groups.

The qualitative agreement between BoSSS and MooNMD also holds for the circularity, see Figure \ref{fig:RBtc2_benchmarkQuantities} in bottom row. For the second test case the benchmark groups agree very well up to $t \approx 1.75$. After that they show large deviations, especially TP2D which is due to the retractions of the filaments after the break up. All other groups do not exhibit a break up of the filaments. 







\section{Conclusion and outlook}
\label{sec:conclusion}

Within this work a two-phase flow solver based on the XDG method is presented, where the approximation space is adapted to be conformal to the position of the interface. An XDG discretization for the transient incompressible Navier-Stokes equations for sharp interface two-phase flows was proposed. For the time discretization an XDG adapted moving interface approach was applied with a third order BDF scheme. Changes in the interface topology and the occurrence of small cut-cells are handled via a cell agglomeration technique. 

The level-set function representing the interface is discretized by two different standard DG fields. One field is used for the discretization and the other for the level-set, resp. interface, evolution. For the latter an extension velocity field is constructed for what a two-staged marching algorithm was presented utilizing an elliptic PDE-based approach on the cut-cells and a fast-marching procedure on the near field cells. The other level-set ensuring the continuity of the interface is formulated by an $L^2$-projection with suitable continuity constraints. 

The presented XDG-solver was validated against a wide range of typical two-phase surface tension driven flow phenomena. For the capillary wave test case the solver showed very good agreement to the analytical solution for a wide range of wave behavior. For higher frequency oscillations the numerical results slightly underestimate the damping rate. Higher order convergence rates could not be observed for the amplitude height in time.

Regarding the stability of the discretization, especially the surface tension force discretization via the Laplace-Beltrami formulation without regularization, transient simulations with a droplet in equilibrium and non-equilibrium state were performed. All simulations showed stable results under mesh refinement studies. However, considering the equilibrium state and evaluating the surface divergence indicates some stability issues for long simulations times.

At last the XDG-solver was tested against other numerical codes within the rising bubble benchmark and showed good to very good agreement to the
benchmark groups for both settings considering the terminal bubble shape and scalar measure quantities, such as the center of mass, rise velocity and circularity. The largest deviation from the benchmark groups was recognized for the circularity within the first setting.

All presented simulations were solved by a direct solver (MUMPS), but considering the simulation of 3D-problems the number of degrees of freedom requires the need of iterative solvers. XDG adapted Newton-methods for the
non-linear Navier-Stokes problem coupled with the level-set evolution are subject of ongoing work at the department of fluid dynamics. However, one should note that the presented method is readily extendable for the three dimensional case. 

\section*{Acknowledgments}
The work of M. Smuda was to some part funded by the 'Excellence Initiative' of the German Federal and State Governments and the Graduate School of Computational Engineering at Technical University Darmstadt and to another part from the German Research Foundation (Deutsche Forschungsgemeinschaft, DFG) within the Collaborative Research Center 1194 ``Interaction between Transport and Wetting Processes'' - Project-ID 265191195.

\appendix

\section{Discretization of the free slip boundary condition}
\label{app:FreeSlipBC}

For the simulation of the rising bubble benchmark in Section \ref{sec:RB} one needs to enforce the free-slip boundary condition \eqref{eq:freeslip}. In context of our discretization \eqref{eq:variationalForm_NSE} the SIP form $a(\veloc, \testV)$ reduces to the components in the normal direction at the free-slip wall $\slip{\edge}$. Thus the bilinear form $a(\veloc, \testV)$ extends to $\tilde a(\veloc, \testV)$ with
\begin{equation}
\begin{aligned}
	\tilde a(\veloc, \testV) = a(\veloc, \testV) + \int_{\slip{\edge}} \visc \left( \normalGam \cdot \left( \gradH{\veloc} + \gradHT{\veloc} \right) \normalGam \right) \cdot 	(\testV \cdot \normalGam) \dif{S}\\
	+ \int_{\slip{\edge}} \visc \left( \normalGam \cdot \left( \gradH{\testV} + \gradHT{\testV} \right) \normalGam \right) \cdot (\veloc \cdot \normalGam) \dif{S}
	- \oint_{\slip{\edge}} \eta (\veloc \cdot \normalGam) \cdot (\testV \cdot \normalGam) \dif{S}.
	\label{eq:viscousTermsCL}
\end{aligned}
\end{equation}
For all other linear forms the free-slip condition corresponds to the Dirichlet boundary condition $\diri{\edge}$.

\section{Additional results for the rising bubble benchmark}
\label{app:risingBubble}

In Hysing et al.\cite{hysing_quantitative_2009} two more error norms are considered besides the $l_2$-error norm \eqref{eq:ScalarErrorNorms2}:
\begin{subequations}
\begin{align}
	l_1 \textrm{-error}: \norm{e}_1 &= \frac{\sum_{n=1}^{\textrm{NTS}} \abs{q_n - q_{n,\textrm{ref}}} }{\sum_{n=1}^{\textrm{NTS}} \abs{q_{n,\textrm{ref}}}}, \label{eq:ScalarErrorNorms1}\\ 
	l_\infty \textrm{-error}: \norm{e}_\infty &= \frac{\max_n \abs{q_n - q_{n,\textrm{ref}}} }{\max_n \abs{q_{n,\textrm{ref}}}}. \label{eq:ScalarErrorNormsInf}
\end{align}
\label{eq:ScalarErrorNorms}
\end{subequations}

\begin{figure}
	\centering
	\begin{subfigure}{.65\textwidth}	
        \includegraphics[width=.9\linewidth]{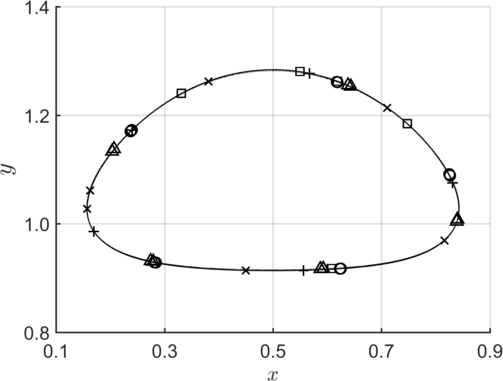}
    \end{subfigure}	
    \begin{subfigure}{.15\textwidth}	
        \includegraphics[width=.9\linewidth]{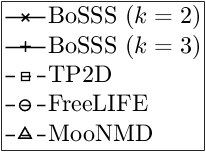}
    \end{subfigure}	
	\caption{Comparison of the terminal bubble shapes at $t=3$ for the test case 1. The resolution of the presented solutions are: BoSSS ($\Pdeg = 2$) with $1/\grdSz = 80$, BoSSS ($\Pdeg = 3$) with $1/\grdSz = 60$, TP2D with $1/\grdSz = 320$, FreeLIFE with $1/\grdSz = 160$, MooNMD with $\textrm{NDOF}_\textrm{int} = 900$.}
	\label{app:RBtc1_terminalShape}
\end{figure}

\begin{center}
\begin{table}
\caption{$l_1$-error norms and $\textrm{ROC}$ for the mesh study of the rising bubble benchmark test case 1.}
\centering
	\begin{tabular}{cc|cc|cc|cc}
		&& \multicolumn{2}{c}{Center of mass} & \multicolumn{2}{c}{Rise velocity} & \multicolumn{2}{c}{Circularity} \\
		&$1/\grdSz$ & $\norm{e}_1$ & $\textrm{ROC}$ & $\norm{e}_1$ & $\textrm{ROC}$ & $\norm{e}_1$ & $\textrm{ROC}$ \\
		\hline
		\hline 
		\multirow{3}{*}{\rotatebox[origin=c]{90}{$\Pdeg = 2$}} & 20 & $4.05 \cdot 10^{-3}$ & - & $9.11 \cdot 10^{-3}$ & - & $4.95 \cdot 10^{-4}$ & - \\
		& 40 & $1.04 \cdot 10^{-3}$ & 1.96 & $2.69 \cdot 10^{-3}$ & 1.76 & $4.07 \cdot 10^{-4}$ & 0.28 \\
		& 60 & $3.73 \cdot 10^{-4}$ & 2.53 & $9.51 \cdot 10^{-4}$ & 2.56 & $1.40 \cdot 10^{-4}$ & 2.63 \\
		\hline 
		\multirow{3}{*}{\rotatebox[origin=c]{90}{$\Pdeg = 3$}} & 10 & $1.21 \cdot 10^{-2}$ & - & $1.37 \cdot 10^{-2}$ & - & $7.29 \cdot 10^{-4}$ & -  \\
		& 20 & $1.78 \cdot 10^{-3}$ & 2.76 & $5.41 \cdot 10^{-3}$ & 1.34 & $1.25 \cdot 10^{-3}$ & -0.78 \\
		& 40 & $5.59 \cdot 10^{-4}$ & 1.67 & $1.53 \cdot 10^{-3}$ & 1.82 & $3.45 \cdot 10^{-4}$ & 1.86
		\label{app:RB_convStudy_errorNorms1}
	\end{tabular}
\end{table}
\end{center}

\begin{center}
\begin{table}
\caption{$l_\infty$-error norms and $\textrm{ROC}$ for the mesh study of the rising bubble benchmark test case 1.}
\centering
	\begin{tabular}{cc|cc|cc|cc}
		&& \multicolumn{2}{c}{Center of mass} & \multicolumn{2}{c}{Rise velocity} & \multicolumn{2}{c}{Circularity} \\
		&$1/\grdSz$ & $\norm{e}_\infty$ & $\textrm{ROC}$ & $\norm{e}_\infty$ & $\textrm{ROC}$ & $\norm{e}_\infty$ & $\textrm{ROC}$ \\
		\hline
		\hline 
		\multirow{3}{*}{\rotatebox[origin=c]{90}{$\Pdeg = 2$}} & 20 & $6.67 \cdot 10^{-3}$ & - & $1.43 \cdot 10^{-2}$ & - & $2.61 \cdot 10^{-3}$ & - \\
		& 40 & $1.74 \cdot 10^{-3}$ & 1.94 & $3.50 \cdot 10^{-3}$ & 2.03 & $8.57 \cdot 10^{-4}$ & 1.61 \\
		& 60 & $6.27 \cdot 10^{-4}$ & 2.52 & $1.20 \cdot 10^{-3}$ & 2.64 & $3.12 \cdot 10^{-4}$ & 2.49 \\
		\hline 
		\multirow{3}{*}{\rotatebox[origin=c]{90}{$\Pdeg = 3$}} & 10 & $2.05 \cdot 10^{-2}$ & - & $2.06 \cdot 10^{-2}$ & - & $2.04 \cdot 10^{-3}$ & -  \\
		& 20 & $2.95 \cdot 10^{-3}$ & 2.80 & $7.38 \cdot 10^{-3}$ & 1.48 & $2.59 \cdot 10^{-3}$ & -0.34 \\
		& 40 & $9.17 \cdot 10^{-4}$ & 1.69 & $2.07 \cdot 10^{-3}$ & 1.83 & $7.71 \cdot 10^{-4}$ & 1.40
		\label{app:RB_convStudy_errorNormsInfty}
	\end{tabular}
\end{table}
\end{center}

\begin{center}
\begin{table}
\caption{Benchmark quantities (test case 1) at distinct values in time for the finest solutions of the respective groups: BoSSS for $\Pdeg = \{2,3\}$, Heimann et al.\cite{heimann_unfitted_2013} for the space pair $(\vectr{X}^{k,2},M^{k,1})$ and the benchmark groups TP2D, FreeLIFE and MooNMD. The values denote: the minimum circularity $\cancel{c}_\textrm{min}$ and the corresponding point in time $t\vert_{\cancel{c} = \cancel{c}_\textrm{min}}$, the maximum rise velocity $V_{c,\textrm{max}}$ and corresponding time $t\vert_{V_c = V_{c,\textrm{max}}}$, and the terminal rise height $y_c(t=3)$.}
\centering
	\begin{tabular}{c|cccccc}
	    & \multicolumn{2}{c}{BoSSS} & Heimann et al.\cite{heimann_unfitted_2013} & \multicolumn{3}{c}{benchmark groups\cite{hysing_quantitative_2009}} \\
		& $\Pdeg = 2$ & $\Pdeg = 3$  & $(\vectr{X}^{k,2},M^{k,1})$ & TP2D & FreeLIFE & MooNMD \\
		& $1/\grdSz = 80$ & $1/\grdSz = 60$  &  $1/\grdSz = 80$  & $1/\grdSz = 320$ & $1/\grdSz = 160$ & $\textrm{NDOF}_{\textrm{int}} = 900$\\
		\hline 
		$\cancel{c}_\textrm{min}$ & 0.9021 & 0.9025 & 0.9010 & 0.9013 & 0.9011 & 0.9013\\
		$t\vert_{\cancel{c} = \cancel{c}_\textrm{min}}$ & 1.8900 & 1.8870 & 1.9000 & 1.9041 & 1.8750 & 1.9000\\
		$V_{c,\textrm{max}}$ & 0.2413 & 0.2412 & 0.2415 & 0.2417 & 0.2421 & 0.2417\\
		$t\vert_{V_c = V_{c,\textrm{max}}}$ & 0.9180 & 0.9180 & 0.9211 & 0.9213 & 0.9313 & 0.9239\\
		$y_c(t=3)$ & 1.0800 & 1.0800 & 1.0809 & 1.0813 & 1.0799 & 1.0817
		\label{app:RBtc1_distinctValues}
	\end{tabular}
\end{table}
\end{center}

\begin{figure}
	\centering
	\begin{subfigure}{.4\textwidth}	
        \includegraphics[width=.9\linewidth]{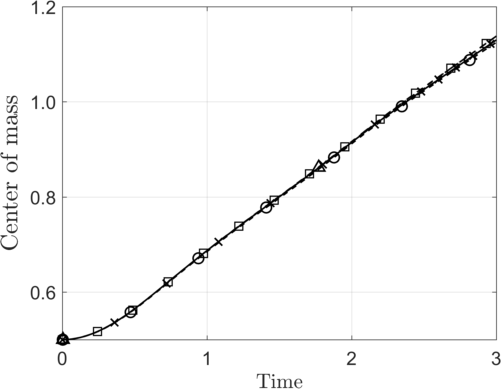}
    \end{subfigure}	
	\begin{subfigure}{.4\textwidth}	
        \includegraphics[width=.9\linewidth]{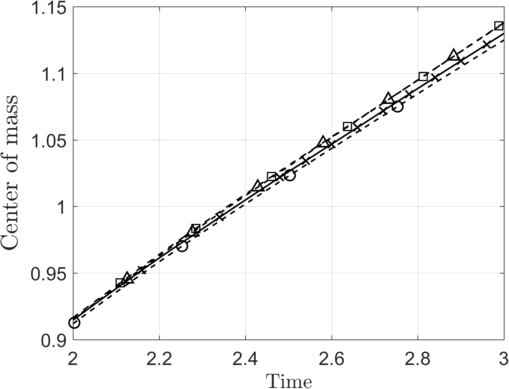}
    \end{subfigure}	
    \begin{subfigure}{.15\textwidth}	
        \includegraphics[width=.9\linewidth]{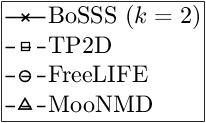}
    \end{subfigure}	
	\caption{Comparison to benchmark groups for the temporal evolution of the $y$-component of the centre of mass (test case 2). On the right the overall evolution is depicted and on the left a zoom for $t = [2, 3]$ is provided. The plotted solutions show the finest solution of each group.}
	\label{fig:RBtc2_CenterOfMass}
\end{figure}

\begin{center}
\begin{table}
\caption{Benchmark quantities (test case 1) at distinct values in time for the finest solutions of the respective groups: BoSSS for $\Pdeg = 2$, Heimann et al.\cite{heimann_unfitted_2013} for the space pair $(\vectr{X}^{k,2},M^{k,1})$ and the benchmark groups TP2D, FreeLIFE and MooNMD. The values denote: the minimum circularity $\cancel{c}_\textrm{min}$ and the corresponding point in time $t\vert_{\cancel{c} = \cancel{c}_\textrm{min}}$, both maximum values for the rise velocity $V_{c,\textrm{max1}}$, $V_{c,\textrm{max2}}$ and the corresponding times $t\vert_{V_c = V_{c,\textrm{max1}}}$, $t\vert_{V_c = V_{c,\textrm{max2}}}$, and the terminal rise height $y_c(t=3)$.}
\centering
	\begin{tabular}{c|ccccc}
	    & BoSSS & Heimann et al.\cite{heimann_unfitted_2013} & \multicolumn{3}{c}{benchmark groups\cite{hysing_quantitative_2009}} \\
		& $\Pdeg = 2$ & $(\vectr{X}^{k,2},M^{k,1})$ & TP2D & FreeLIFE & MooNMD \\
		& $1/\grdSz = 80$ + AMR 1 &  $1/\grdSz = 80$ & $1/\grdSz = 640$ & $1/\grdSz = 160$ & $\textrm{NDOF}_{\textrm{int}} = 900$ \\
		\hline
		$\cancel{c}_\textrm{min}$ & 0.5093 & 0.4903 & 0.5869 & 0.4647 & 0.5144 \\
		$t\vert_{\cancel{c} = \cancel{c}_\textrm{min}}$ & 2.9850 & 3.0000 & 2.4004 & 3.0000 & 3.0000 \\
		$V_{c,\textrm{max}1}$ & 0.2500 & 0.2500 & 0.2524 & 0.2514 & 0.2502 \\
		$t\vert_{V_c = V_{c,\textrm{max}1}}$ & 0.7290 & 0.7281 & 0.7332 & 0.7281 & 0.7317 \\
		$V_{c,\textrm{max}2}$ & 0.2373 & 0.2377 & 0.2434 & 0.2440 & 0.2393 \\
		$t\vert_{V_c = V_{c,\textrm{max}2}}$ & 2.0200 & 2.0281 & 2.0705 & 1.9844 & 2.0600 \\
		$y_c(t=3)$ & 1.1300 & 1.1263 & 1.1380 & 1.1249 & 1.1376
		\label{app:RBtc2_distinctValues}
	\end{tabular}
\end{table}
\end{center}

\bibliographystyle{unsrt}
\bibliography{XNSEpaper}%


\end{document}